\documentclass[11pt,amssymb]{article}
\usepackage{amsfonts,latexsym}

\newcommand{\R}{\mathbb R}

\newcommand{\N}{\mathbb N}

\newcommand{\C}{\mathbb C}

\newtheorem{theorem}{Theorem}
\newtheorem{lemma}{Lemma}

\newtheorem{corollary}{Corollary}
\newtheorem{definition}{Definition}

\def\ep{\varepsilon}

\begin{document}
\title{Bounded solutions for an ordinary differential system from the Ginzburg-Landau theory.}
\author{Anne Beaulieu, LAMA, Univ Paris Est Creteil, Univ Gustave Eiffel, UPEM, CNRS, F-94010, Cr\'eteil, France.}
  \maketitle
 anne.beaulieu@u-pec.fr                    
                                
   {\bf Abstract.} \\
In this paper, we look at a linear system of ordinary differential equations as derived from the two-dimensional Ginzburg-Landau equation. In two cases, it is known that this system admits bounded solutions coming from the invariance of the Ginzburg-Landau equation by translations and rotations. The specific contribution of our work is to prove that in the other cases, the system does not admit any bounded solutions. We show that this bounded solution problem is related to an eigenvalue problem.\\

                          \section{Introduction.}
Let $n$ and $d$ be given integers, $n\geq1$, $d\geq1$. We define the following system
\begin{equation}\label{eq:GL}
\left\{\begin{array}{rl}
a''+{a'\over r}-{(n-d)^2\over r^2}a-f_d^2b&=-(1-2f_d^2)a\\
b''+{b'\over r}-{(n+d)^2\over r^2}b-f_d^2a&=-(1-2f_d^2)b
\end{array}
\right.
\end{equation}
 and the following equations

\begin{equation}\label{eq:GL0}
a''+{a'\over r}-{d^2\over r^2}a=-(1-f_d^2)a
\end{equation}
and

\begin{equation}\label{eq:GLR}
a''+{a'\over r}-{d^2\over r^2}a-{2a}f^2_d=-(1-f_d^2)a.
\end{equation}
with the variable $r>0$, and for real valued functions  $r\mapsto a(r)$  and  $r\mapsto b(r)$.\\
Here $f_d$ is the only solution of the differential equation
\begin{equation}\label{eq:f}
f_d''+{f_d'\over r}-{d^2\over r^2}f_d=-f_d(1-f_d^2)
\end{equation}
with the conditions $f_d(0)=0$ and $\lim_{+\infty}f_d=1$.

 Let us consider the Ginzburg-Landau equation on a bounded connected domain $\Omega$, 
 
\begin{equation}\label{eq:GLep}
\left\{\begin{array}{c}
-\Delta u={1\over\ep^2}u(1-\vert u\vert^2 )\hbox{   in     }\Omega\\
u=g\hbox{    in   }\partial\Omega
\end{array}
\right.
\end{equation}
where $\ep>0$ is a small parameter, $u$ and $g$ have complex values and degree $(g,\partial\Omega)\geq1$. Let us consider the following equation
\begin{equation}\label{eq:GLespace}
\begin{array}{c}
-\Delta u=u(1-\vert u\vert^2)\hbox{   in   }\R^2\\
\end{array}
\end{equation}
where $u$ is a complex valued map. 
The study of the energy-minimizing solutions of equation (\ref{eq:GLep}) is in the book of Bethuel, Brezis H\'elein, \cite{BBH}.\\
Let us explain how the system (\ref{eq:GL}) and the equations (\ref{eq:GL0}) and (\ref{eq:GLR}) are derived from the equations (\ref{eq:GLep}) and (\ref{eq:GLespace}).\\
We denote $\mathcal{N}_{\ep}(u)=\Delta u+{1\over\ep^2}u(1-\vert u\vert^2 )$. Let
  $u_0(x)=f_d({\mid x\mid\over\ep})e^{id\theta}$. We have $\mathcal{N}_{\ep}(u_0)=0$. We will always 
  denote $$f(r)=f_d({r\over\ep}).$$
   We differentiate $\mathcal{N}_{\ep}$ at $u_0$.
$$d\mathcal{N}_{\ep}(u_0)(\omega)=\Delta \omega+{\omega\over\ep^2}(1-f^2)-{2\over\ep^2}f^2 e^{id\theta}e^{id\theta}.\omega,$$
where $\omega$ is any complex valued function and $2u.\omega=\overline{u} \omega+\overline{\omega} u.$ We will use the operator
$e^{-id\theta}d\mathcal{N}_{\ep}(u_0)e^{id\theta}$ instead of $d\mathcal{N}_{\ep}(u_0)$.
We consider the Fourier expansion
$$\omega(x)=\sum_{n\geq1}(a_n(r)e^{-in\theta}+b_n(r)e^{in\theta})+a_0(r),\quad a_n(r)\in \C,\quad b_n(r)\in\C.$$
Letting 
$\omega_n(x)=a_n(r)e^{-in\theta}+b_n(r)e^{in\theta},$
we have
$$2e^{id\theta}.e^{id\theta}\omega_n=\omega_n+\overline{\omega}_n=(b_n+\overline{a}_n)e^{in\theta}+(\overline{b}_n+a_n)e^{-in\theta}.$$
$$\hbox{Moreover} \qquad e^{-id\theta}\Delta(e^{id\theta}\omega)=\Delta \omega-{d^2\over r^2}\omega+i{2d\over r^2}{\partial\omega\over\partial\theta}.$$
Then
$$e^{-id\theta}d\mathcal{N}_{\ep}(u_0)e^{id\theta}\omega=\sum_{n\geq1}e^{-in\theta}(a''_n+{a'_n\over r}-{(n-d)^2\over r^2}a_n+{a_n\over\ep^2}(1-f^2)-{a_n\over\ep^2}f^2-{\overline{b}_n\over\ep^2}f^2)
$$
$$+\sum_{n\geq1}e^{in\theta}( b''_n+{b'_n\over r}-{(n+d)^2\over r^2}b_n+{b_n\over\ep^2}(1-f^2)-{b_n\over\ep^2}f^2-{\overline{a}_n\over\ep^2}f^2)
$$
$$
+ a''_0+{a'_0\over r}-{d^2\over r^2}a_0+{a_0\over\ep^2}(1-f^2)-{a_0+\overline{a}_0\over\ep^2}f^2.$$
Separating the Fourier components of  $e^{-id\theta}d\mathcal{N}_{\ep}(u_0)e^{id\theta}\omega$, we can consider the operators
$$ \hbox{for $n\geq1$},\quad\mathcal{L}_n(a_n,b_n)=\left\{\begin{array}{c}
a''_n+{a'_n\over r}-{(n-d)^2\over r^2}a_n+{a_n\over\ep^2}(1-2f^2)-{\overline{b}_n\over\ep^2}f^2\\
 b''_n+{b'_n\over r}-{(n+d)^2\over r^2}b_n+{b_n\over\ep^2}(1-2f^2)-{\overline{a}_n\over\ep^2}f^2\\
\end{array}
\right.$$
$$\hbox{and, for $n=0$},\quad
\mathcal{L}_0(a_0)=
a''_0+{a'_0\over r}-{d^2\over r^2}a_0+{a_0\over\ep^2}(1-f^2)-{a_0+\overline{a}_0\over\ep^2}f^2.$$
When we have to solve the system  $(\mathcal{L}_n(a_n,b_n)=(\alpha_n,\beta_n),\mathcal{L}_0(a_0)=\alpha_0)$, for some given $(\alpha_n,\beta_n)\in\C\times\C$ and $\alpha_0\in\C$, we are led to consider separately the real part and the imaginary part. So we consider the following operators, where $a_n$ and $b_n$ are real valued functions 
$$\hbox{for $n\geq1$}\quad{\mathcal{L}}_{n,\mathcal{R}}:(a_n,b_n)\mapsto\left\{\begin{array}{c}
a''_n+{a'_n\over r}-{(n-d)^2\over r^2}a_n+{a_n\over\ep^2}(1-2f^2)-{{b}_n\over\ep^2}f^2\\
 b''_n+{b'_n\over r}-{(n+d)^2\over r^2}b_n+{b_n\over\ep^2}(1-2f^2)-{{a}_n\over\ep^2}f^2\\
\end{array}
\right.;$$
$${\mathcal{L}}_{n,\mathcal{I}}:(a_n,b_n)\mapsto\left\{\begin{array}{c}
a''_n+{a'_n\over r}-{(n-d)^2\over r^2}a_n+{a_n\over\ep^2}(1-2f^2)+{{b}_n\over\ep^2}f^2\\
 b''_n+{b'_n\over r}-{(n+d)^2\over r^2}b_n+{b_n\over\ep^2}(1-2f^2)+{{a}_n\over\ep^2}f^2\\
\end{array}
\right.$$
$$\hbox{and, for $n=0$},\quad
{\mathcal{L}}_{0,\mathcal{I}}:a_0\mapsto
a''_0+{a'_0\over r}-{d^2\over r^2}a_0+{a_0\over\ep^2}(1-f^2)\quad;$$
$$\quad{\mathcal{L}}_{0,\mathcal{R}}: a_0\mapsto
a''_0+{a'_0\over r}-{d^2\over r^2}a_0+{a_0\over\ep^2}(1-f^2)-{2a_0\over\ep^2}f^2.$$
Considering $(-a_n,b_n)$, we see that, for $n\geq1$, only one of the operators ${\mathcal{L}}_{n,\mathcal{R}}$ or ${\mathcal{L}}_{n,\mathcal{I}}$ is of interest. 
The eigenvalue problems ${\mathcal{L}}_{n,\mathcal{R}}(a_n,b_n)=-\lambda(\ep)(a_n, b_n),$ $(a_n,b_n)(1)=0$, for all integers $d\geq1$ and $n\geq1$, as well as the problems ${\mathcal{L}}_{0,\mathcal{R}}(a_0)=-\lambda(\ep)a_0$ and ${\mathcal{L}}_{0,\mathcal{I}}(a_0)=-\lambda(\ep)a_0$, $a_0(1)=0$, have been studied in several papers, including \cite{LiebLoss}, \cite{Miro}, \cite{TCL1}, \cite{TCL2}, \cite{SR}. In the third chapter of their book \cite{PR}, Pacard and Rivi\`ere study the system (\ref{eq:GL}) and the equations (\ref{eq:GL0}) and (\ref{eq:GLR}) for $d=1$. These authors' aim is to construct some solutions for (\ref{eq:GLep}).\\
Let us now bring together some of the results contained in the above studies.
\begin{theorem}\label{essentiel}
For $d\geq1$ and $n\geq1$, the existence of an eigenvalue $\lambda(\ep)\rightarrow0$ as $\ep\rightarrow0$ is equivalent to the existence of a bounded solution of (\ref{eq:GL}).
For the equations (\ref{eq:GL0}) and (\ref{eq:GLR}) and for all $d\geq1$, the results of \cite{PR} are valid for all $d\geq1$, that is the real vector space of the bounded solutions of (\ref{eq:GL0}) is one-dimensional, spanned by $f_d$ and there is no bounded solution of (\ref{eq:GLR}). For $n=1$, the vector space of the bounded solutions of (\ref{eq:GL}) is also a one dimensional vector space, spanned by 
$(f_d'+{d\over r}f_d, f_d'-{d\over r}f_d)$. For $d=1$ and $n\geq2$, there are no bounded solutions. For $d\geq2$ and for $n\geq 2d-1$, there are no 
bounded solutions.
\end{theorem}
A bounded solution is any solution defined at $r=0$ and which has a finite limit as $r\rightarrow+\infty$. For  all $d\geq1$, the known bounded solutions, for $n=0$ and $n=1$, come from the invariance of the Ginzburg-Landau equation with respect to the translations and the rotations.\\
The present paper's aim is to prove the following
\begin{theorem}\label{pasdesolutionbornee}
For all real numbers $d$ and $n$ such that $d\geq1$ and $n>1$, the system (\ref{eq:GL}) has no bounded solution.
\end{theorem}
We will consider $n$ and $d$ as real parameters, although the Ginzburg-Landau problem is about integers $n$ and $d$. So we have to consider the functions $f_d$, for $d\in\R$, $d\geq1$. But in \cite{HH}, where all the solutions of (\ref{eq:f}) are studied, the authors consider only the case $d\in\N^{\star}$. However, this hypothesis is not essential in their paper. Here are the properties of $f_d$ we need.
\begin{theorem}\label{fd}
Let $d\in\R^{+\star}$. For all $a>0$ there exists a unique solution of (\ref{eq:f}) such that $\lim_{r\rightarrow0}{1\over a}f(r)r^{-d} =1$.
There exists a unique value $A_d>0$ such that this solution is defined in $\R^{+\star}$ and non decreasing. 
Denoting it by $f_d$, we have the expansions
\begin{equation}\label{eq:fenlinfini}
f_d(r)=1-{d^2\over 2r^2}+O({1\over r^4})\hbox{   near $+\infty$}
\end{equation}
and
\begin{equation}\label{eq:fen0}
f_d(r)=A_d(r^d-{1\over 4(d+1)}r^{d+2})+O(r^{d+4})\hbox{   near 0}.
\end{equation}
Moreover, if we denote $g^{(d)}(r)=r^{-d}f_d(r)$, then, for all $\alpha>0$, the map $d\mapsto g^{(d)}$ is continuous from $]0,+\infty[$ into $L^{\infty}([0,\alpha])$.
\end{theorem}
For $d\in\R^{+\star}$, $\gamma_1\in\R^{+}$ and $\gamma_2\in\R^{+}$, we define the following system
\begin{equation}\label{eq:1}
\left\{\begin{array}{rl}
a''+{a'\over r}-{\gamma_1^2\over r^2}a-f_d^2b-f_d^2a&=-(1-f_d^2)a\\
b''+{b'\over r}-{\gamma_2^2\over r^2}b-f_d^2a-f_d^2b&=-(1-f_d^2)b.
\end{array}
\right.
\end{equation}
Letting $x=a+b$ and $y=a-b$, we are led to the system verified by $(x,y)$, that is
\begin{equation}\label{eq:2}
\left\{\begin{array}{ll}
x''+{x'\over r}-{\gamma^2\over r^2}x+{\xi^2\over r^2}y-2f_d^2x&=-(1-f_d^2)x\\
y''+{y'\over r}-{\gamma^2\over r^2}y+{\xi^2\over r^2}x&=-(1-f_d^2)y
\end{array}
\right.
\end{equation}
with 
$$\gamma^2={\gamma_1^2+\gamma_2^2\over 2}\quad\hbox{and}\quad \xi^2={\gamma_2^2-\gamma_1^2\over2}.$$
 In all the paper, when there is no other indication, we will suppose that $d\geq1$ and ${\gamma^2_1+\gamma^2_2\over 2}-d^2\geq1$. We will denote
\begin{equation}\label{eq:defn}
n=\sqrt{{\gamma^2_1+\gamma^2_2\over 2}-d^2}.
\end{equation}
But Theorem \ref{solutionscontinuesen0}, Theorem \ref{solutions continues en infty} and Theorem \ref{lesextremes} will be valid for $d>0$ and ${\gamma^2_1+\gamma^2_2\over 2}-d^2>0$. \\
Let us cite a supposedly well-known principle : for any real number $R>0$ and for any given Cauchy data $(a(R),a'(R), b(R), b'(R))$, the system (\ref{eq:1}) has a unique solution, defined in $\R^{+\star}$. This solution is continuous wrt the real positive parameters $d$, $\gamma_1$ and $\gamma_2$, because the coefficients of the system depend continuously on them. Moreover, when the Cauchy data depends continuously on the parameters, so does the solution $(a(r),a'(r), b(r), b'(r))$, which, consequently, is bounded independently of $r$ and of the parameters, when $r$ and the parameters stay in a given compact set. This principle comes from the Cauchy-Lipschitz Theorem, whose proof rests on an application of the Banach Fixed Point Theorem to a suitable integral equation. However, we don't know whether 
a given solution keeps the same behavior at 0 or at $+\infty$ for all the values of the parameters, even when this solution is continuous wrt the parameters. So we begin with the definition of some continuous solutions wrt to the parameters in a certain range, and whose behaviors remain inchanged, either at 0 or at $+\infty$. \\
To begin with, let us give the following definition
\begin{definition}\label{Oeto}
We say that
\begin{enumerate}
\item
$a=$O$(f)$  at 0
if there exists  $R>0$ and $C>0$ such that
$$\forall r\in]0,R],\quad \vert a(r)\vert\leq C\vert f(r)\vert.$$
\item
$a$  has the behavior $f$ at 0, and we denote $a\sim_0f$,
if there exists a map $g$, such that
$$\lim_0 g=0,\quad \vert a-f\vert=\hbox{O}(fg).$$
\item
$a=o(f)$ at 0 if there exists a map $g$, such that
$$\lim_0 g=0,\quad a=fg.$$
\end{enumerate}
We will use the same convention at $+\infty$.
 \end{definition}
We will consider that $(d,\gamma_1,\gamma_2)$ belongs to the set
$$\mathcal{D}=\{(d,\gamma_1,\gamma_2)\in(\R_+)^3;d\geq1;\gamma_2>1; \quad0\leq\gamma_1\leq\gamma_2<\gamma_1+2d+2\}.$$
Let us remark that $(d,\vert n-d\vert,n+d)\in\mathcal{D}$, whenever $d\geq 1$ and $n\geq1$, $d\in\R$, $n\in\R$. \\
We will need the following subsets of $\mathcal{D}$
$$\mathcal{D}_1=\{(d,\gamma_1,\gamma_2)\in\mathcal{D};\gamma_1>0\}$$
and
\begin{equation}\label{eq:defD}
\mathcal{D}_2=\{(d,\gamma_1,\gamma_2)\in\mathcal{D};0\leq\gamma_1<{1\over4}; -\gamma_1- \gamma_2+2d+2>0;-\gamma_2+2d+1>0\}.
 \end{equation}
  Whenever $d\geq 1$ and $n\geq1$, $n\in\R$, $d\in\R$, we remark that $(d,\vert n-d\vert,n+d)\in\mathcal{D}_1$ when $n\neq d$ and that $(d,\vert n-d\vert,n+d)\in\mathcal{D}_2$ when 
 $\vert n-d\vert<{1\over 4}$.\\
 The following theorem is about a base of solutions defined near 0.
 \begin{theorem}\label{solutionscontinuesen0}
For all $(d,\gamma_1, \gamma_2)\in\mathcal{D}$,
there exist four independent solutions $(a,b)$ of (\ref{eq:1}) verifying the following conditions
\begin{enumerate}
\item 
$(a_1(r),b_1(r))\sim_0 (O(r^{\gamma_2+2d+2}),r^{\gamma_2})\hbox{  and  }(a'_1(r),b'_1(r))\sim_0 (O(r^{\gamma_2+2d+1}),\gamma_2r^{\gamma_2-1}).$
\item
$(a_2(r),b_2(r))\sim_0\left\{\begin{array}{c}(O(r^2\theta(r)),r^{-\gamma_2})\quad\hbox{  if   }\quad(d,\gamma_1,\gamma_2)\in\mathcal{D}_1\\
(O(r^{-\gamma_2+2d+2}),r^{-\gamma_2})\quad\hbox{  if   }\quad(d,\gamma_1,\gamma_2)\in\mathcal{D}_2\\
\end{array}
\right.
$\\
$(a'_2(r),b'_2(r))\sim_0\left\{\begin{array}{c}(O(r\theta(r)),-\gamma_2r^{-\gamma_2-1})\quad\hbox{  if   }\quad(d,\gamma_1,\gamma_2)\in\mathcal{D}_1\\
(O(r^{-\gamma_2+2d+1}),-\gamma_2r^{-\gamma_2-1})\quad\hbox{  if   }\quad(d,\gamma_1,\gamma_2)\in\mathcal{D}_2\\
\end{array}
\right.
$\\
$\hbox{where  }\theta(r)=\left\{\begin{array}{c}{-r^{\gamma_1-2}+r^{-\gamma_2+2d}\over \gamma_1+\gamma_2-2d-2}\quad\hbox{  if  }\quad\gamma_1+\gamma_2-2d-2\neq0\\
-r^{\gamma_1-2}\log r\quad\hbox{  if  }\quad\gamma_1+\gamma_2-2d-2=0.
\end{array}
\right.
$
\item
$(a_3(r),b_3(r))\sim_0(r^{\gamma_1},O(r^{\gamma_1+2d+2}))\hbox{  and, if $\gamma_1\neq0$  }    (a'_3(r),b'_3(r))\sim_0(\gamma_1r^{\gamma_1-1},O(r^{\gamma_1+2d+1}))$\\
while, if $\gamma_1=0$,  $(a'_3(r),b'_3(r))=(O(r),O(r^{2d+1}))$.
\item 
$(a_4(r),b_4(r))\sim_0\left\{\begin{array}{c}
(r^{-\gamma_1}, O(r^2\tilde{\theta}(r))\quad\hbox{  if   }\quad(d,\gamma_1,\gamma_2)\in\mathcal{D}_1\\
(\tau(r),O(\tau(r)r^{2d+2}))\quad\hbox{  if   }\quad(d,\gamma_1,\gamma_2)\in\mathcal{D}_2
\end{array}
\right.
$\\
and
$(a'_4(r),b'_4(r))\sim_0\left\{\begin{array}{c}
(r^{-\gamma_1-1}, O(r\tilde{\theta}(r))\quad\hbox{  if   }\quad(d,\gamma_1,\gamma_2)\in\mathcal{D}_1\\
(\tau'(r),O(\tau'(r)r^{2d+2}))\quad\hbox{  if   }\quad(d,\gamma_1,\gamma_2)\in\mathcal{D}_2
\end{array}
\right.
$\\
where\\ $\tilde{\theta}(r)=\left\{\begin{array}{c}{-r^{\gamma_2-2}+r^{-\gamma_1+2d}\over \gamma_1+\gamma_2-2d-2}\hbox{  if  }\gamma_1+\gamma_2-2d-2\neq0\\
-r^{\gamma_2-2}\log r\hbox{  if  }\gamma_1+\gamma_2-2d-2=0
\end{array}
\right.
;\tau(r)=\left\{\begin{array}{c}{r^{-\gamma_1}-r^{\gamma_1}\over 2\gamma_1}\hbox{  if  } \gamma_1\neq0 \\
-\log r\hbox{  if  } \gamma_1=0.
\end{array}
\right.
$
\item 
For $j=1$ and for $j=3$, for all $r>0$, the maps \\
$(d,\gamma_1,\gamma_2)\mapsto(a_j(r),a'_j(r),b_j(r),b'_j(r))\hbox{  are continuous in  }\mathcal{D}.$
\item
For $j=1$ and for $j=3$, and for all $r>0$, $(a_j(r),a'_j(r),b_j(r),b'_j(r))$ is differentiable wrt to $\gamma_1$ and wrt $\gamma_2$, whenever $(d,\gamma_1,\gamma_2)\in\mathcal{D}$, and $\gamma_2>\gamma_1$.\\
Moreover the map $(d,\gamma_1,\gamma_2)\mapsto {\partial\over\partial\gamma_i} (a_j(r),a'_j(r),b_j(r),b'_j(r))$ is continous, for $i=1$ and $i=2$. 
And we have
\begin{equation}\label{eq:cptderivee1}
({\partial a_1\over\partial\gamma_i},{\partial a'_1\over\partial\gamma_i},{\partial b_1\over\partial\gamma_i},{\partial b'_1\over\partial\gamma_i})(r)\sim_0\log r( O(r^{\gamma_2+2d+2}),O(r^{\gamma_2+2d+1}),r^{\gamma_2},\gamma_2r^{\gamma_2-1})
\end{equation}
and, if $\gamma_1\neq0$
\begin{equation}\label{eq:cptderivee3gammanonnul}
({\partial a_3\over\partial\gamma_i},{\partial a'_3\over\partial\gamma_i},{\partial b_3\over\partial\gamma_i},{\partial b'_3\over\partial\gamma_i})(r)
\end{equation}
$\sim_0\log r( r^{\gamma_1},\gamma_1r^{\gamma_1-1}+O(r^{\gamma_1+1 }),O(r^{\gamma_1+2d+2}),O(r^{\gamma_1+2d+1}))$
\item
For $j=2$ or for $j=4$, the same notation $(a_j,b_j)$ is used for two solutions, one of them being defined for  $(d,\gamma_1,\gamma_2)\in\mathcal{D}_1$, the other one being defined for  $(d,\gamma_1,\gamma_2)\in\mathcal{D}_2$. \\
Moreover, for each domain $\mathcal{D}_i$, $i=1,2$ and for all $r>0$
the maps
$(d,\gamma_1,\gamma_2)\mapsto (a_j(r),a'_j(r),b_j(r), b'_j(r))$ are continuous in $\mathcal{D}_i$. For each $r>0$, the partial differentiability of $(a_j(r),a'_j(r),b_j(r), b'_j(r))$ wrt $\gamma_1$ or wrt $\gamma_2$ is also true separatly in each domain $\mathcal{D}_i$, $i=1,2$.\\
\end{enumerate}
\end{theorem}
The second theorem is about a base of solutions defined near $+\infty$.
\begin{theorem}\label{solutions continues en infty}
\begin{enumerate}
\item
 We have a base of four solutions $(a,b)$ of (\ref{eq:1}), with given behaviors at $+\infty$. In order to distinguish these solutions from the solutions defined in Theorem \ref{solutionscontinuesen0}, we use the notation $(u_i,v_i)$, $i=1,\ldots,4,$ for these solutions. We have 
$$(u_1(r),v_1(r))\sim_{r\rightarrow+\infty}({e^{\sqrt 2 r}/\sqrt r},{e^{\sqrt 2 r}/\sqrt r})(1+O(r^{-2}));$$
$$ (u_2(r),v_2(r))\sim_{r\rightarrow+\infty}({e^{-\sqrt 2 r}/\sqrt r},{e^{-\sqrt 2 r}/\sqrt r})(1+O(r^{-2}));$$
and
$$(u_3(r),v_3(r))\sim_{r\rightarrow+\infty}(r^{-n},-r^{-n})(1+O(r^{-2}));$$
$$(u_4(r),v_4(r))\sim_{r\rightarrow+\infty}(r^{n},-r^{n})(1+O(r^{-2})).$$
\item
Except for $j=2$, the construction of $(u_j,v_j)$ is done separatly for each compact subset $ \mathcal{K}$ of $\mathcal{D}$. For each of the four solutions and for all $r>0$ the map  $(d,\gamma_1,\gamma_2)\mapsto (u_j(r),u'_j(r),v_j(r)),v'_j(r))$ is continuous on $\mathcal{K}$. There partial derivatives wrt $\gamma_1$ and wrt $\gamma_2$ exist whenever $\gamma_1<\gamma_2$ and are continuous. We have\\
$({\partial u_1\over\partial \gamma_i},{\partial u'_1\over\partial \gamma_i},{\partial v_1\over\partial \gamma_i}, {\partial v'_1\over\partial \gamma_i})(r)
\sim_{r\rightarrow+\infty}{e^{\sqrt2 r}\over\sqrt r}\log r(O(r^{-2}),O(r^{-3}),O(r^{-2}),O(r^{-3}));
$\\
$({\partial u_2\over\partial \gamma_i},{\partial u'_2\over\partial \gamma_i},{\partial v_2\over\partial \gamma_i}, {\partial v'_2\over\partial \gamma_i})(r)
\sim_{r\rightarrow+\infty}{e^{-\sqrt2 r}\over\sqrt r}\log r(O(r^{-2}),O(r^{-3}),O(r^{-2}),O(r^{-3}))
$;\\
$({\partial u_3\over\partial \gamma_i},{\partial u'_3\over\partial \gamma_i},{\partial v_3\over\partial \gamma_i}, {\partial v'_3\over\partial \gamma_i})(r)
\sim_{r\rightarrow+\infty}\log r(r^n,O(r^{n-1}),-r^n,O(r^{n-1}))(1+O(r^{-2}))
$;\\
$({\partial u_4\over\partial \gamma_i},{\partial u'_4\over\partial \gamma_i},{\partial v_4\over\partial \gamma_i}, {\partial v'_4\over\partial \gamma_i})(r)
\sim_{r\rightarrow+\infty}\log r(r^{-n},O(r^{-n-1}),-r^{-n},O(r^{-n-1}))(1+O(r^{-2})).
$
\end{enumerate}
\end{theorem}
By our construction, the solution $(u_j,v_j)$ depends on the given compact set $\mathcal{K}$, except for $j=2$. 
For $j=1$, this difficulty disappears after the proof of Theorem \ref{lesextremes}. For the other solutions, named $(u_3,v_3)$ and $(u_4,v_4)$, we will have to make sure that the parameter $(d,\gamma_1,\gamma_2)$ stays in a compact set, as soon as we want and use the continuity and the differentiability of these solutions wrt the parameters.\\
In \cite{SR} we already gave the behaviors of a base of solutions at 0 and at $+\infty$.
But, in the present paper, the continuity wrt to $(d,\gamma_1,\gamma_2)$, especially of the five solutions $(a_3,b_3)$ and $(a_1,b_1)$ (defined at 0) and $(u_2,v_2)$, $(u_3,v_3)$, $(u_4,v_4)$ (defined at $+\infty$) and there differentiability wrt $\gamma_1$ and $\gamma_2$, are essential. \\
The following theorem connects the least behavior at 0 to the exponentially blowing up behavior at $+\infty$ and the least behavior at $+\infty$ to the greater blowing up behavior at 0.
\begin{theorem}\label{lesextremes}
Suppose that $d>0$ and that $\gamma_2\geq\gamma_1\geq0$, $(\gamma_2^2+\gamma_1^2)/2> d^2$. Let $(a_1,b_1)$ be the solution of (\ref{eq:GL}) defined by
$(a_1,b_1)\sim_0(O(r^{\gamma_2+2d+2}),r^{\gamma_2})$. Then $(a_1,b_1)$ blows up exponentially at $+\infty$. Let $(u_2,v_2)$ be the solution of (\ref{eq:GL}) defined by  $(u_2,v_2)\sim_{+\infty}({e^{-\sqrt2 r}\over\sqrt r},{e^{-\sqrt2 r}\over\sqrt r})$. Then $(u_2,v_2)\sim_0 C(o(r^{-\gamma_2}),r^{-\gamma_2})$, for some $C\neq0$.
\end{theorem}
Now, let us relate the problem (\ref{eq:1}) to an eigenvalue problem, which is a little bit different from the one considered in the previous works on the subject, but, for our proof, we find it more suitable.\\
Let $0\leq\gamma_1<\gamma_2$, $\mu\in\R$ and $\ep>0$ be given and let us consider the following system
 \begin{equation}\label{eq:GLvp}
\left\{\begin{array}{rl}
a''+{a'\over r}-{\gamma_1^2\over r^2}a-{1\over\ep^2}f^2a-{1\over\ep^2}f^2b&=-{1\over\ep^2}\mu(1-f^2)a\\
b''+{b'\over r}-{\gamma_2^2\over r^2}b-{1\over\ep^2}f^2b-{1\over\ep^2}f^2a&=-{1\over\ep^2}\mu(1-f^2)b
\end{array}
\right.
\end{equation}
for $r\in]0,1]$, with the condition 
$a(1)=b(1)=0.$
Let us explain in which sense this can be considered as an eigenvalue problem.\\
Let ${\gamma_1}\geq0$ be given. We define
$$\mathcal{H}_{\gamma_1}=\{r\mapsto (a(r),b(r))\in\R^2; (ae^{i\gamma_1\theta},be^{i\theta})\in H_0^1(B(0,1),\C)\times H_0^1(B(0,1),\C)\},$$
where $(r,\theta)$ are the polar coordinates in $\R^2$. The dependence on $\gamma_1$ is needed to distinguish $\gamma_1=0$ and $\gamma_1\neq0$.
We endow $\mathcal{H}_{\gamma_1}$ with the scalar product
$$<(a,b)\vert(u,v)>=\int_0^1(ra'u'+rb'v'+{\gamma^2_1\over r}au+{1\over r}bv)dr$$
and then $\mathcal{H}_{\gamma_1}$ is a Hibert space.
Let $\mathcal{H}'_{\gamma_1}$  be the topological dual space of $\mathcal{H}_{\gamma_1}$.\\
We consider the $\mathcal{T}_{\gamma_1,\gamma_2}:\mathcal{H}_{\gamma_1}\rightarrow \mathcal{H}'_{\gamma_1}$
defined by
$$<\mathcal{T}_{\gamma_1,\gamma_2}(a,b),(u,v)>_{\mathcal{H}',\mathcal{H}}
=\int_0^1(ra'u'+rb'v'+{\gamma_1^2\over r}au+{\gamma_2^2\over r}bv+{r\over\ep^2}f^2(a+b)(u+v))dr.$$
We remark that
$$((a,b),(u,v))\mapsto <\mathcal{T}_{\gamma_1,\gamma_2}(a,b),(u,v)>_{\mathcal{H}'_{\gamma_1},\mathcal{H}_{\gamma_1}}$$
is a scalar product on $\mathcal{H}_{\gamma_1}$. So, $\mathcal{T}_{\gamma_1,\gamma_2}$ is an isomorphism, by the Riesz Theorem.\\
Last, let us define the embedding 
$$\begin{array}{rl}I:&\mathcal{H}_{\gamma_1}\rightarrow\mathcal{H}'_{\gamma_1}\\
&(a,b)\mapsto ((u,v)\mapsto\int_0^1r(au+bv)dr)
\end{array}
$$
Since the embedding $H_0^1(B(0,1))\times H_0^1(B(0,1))\subset L^2(B(0,1))\times L^2(B(0,1))$ is compact, then $I$ is compact.\\
Let us define
$
 \mathcal{C}={1\over\ep^2}(1-f^2)I.
$
Since $\mathcal{C}$ is a compact operator and thanks to the continuity of $\mathcal{T}_{\gamma_1,\gamma_2}^{-1}$, then $\mathcal{T}_{\gamma_1,\gamma_2}^{-1}\mathcal{C}$ is a compact operator from $\mathcal{H}_{\gamma_1}$ into itself. By the standard theory of self adjoint compact operators, there exists a Hilbertian base of $\mathcal{H}$ formed of eigenvectors of $\mathcal{T}_{\gamma_1,\gamma_2}^{-1}\mathcal{C}.$\\
Now let us define $m_{\gamma_1,\gamma_2}(\ep)$ as the first eigenvalue for the above eigenvalue problem in $\mathcal{H}_{\gamma_1}$, that is 
\begin{equation}
{m}_{\gamma_1,\gamma_2}(\ep)=\inf_{(a,b)\in \mathcal{H}^2_{\gamma_1}/\{(0,0)\}}{\int_0^1(ra'^2+rb'^2+{\gamma_1^2\over r}a^2+{\gamma_2^2\over r}b^2+{r\over\ep^2}f_d^2({r\over\ep})(a+b)^2)dr\over{1\over\ep^2}\int_0^1r(1-f_d^2({r\over\ep}))(a^2+b^2)dr}
\end{equation}
and let us define 
\begin{equation}\label{eq:defm0}
m_0(\ep)=\inf_{a\in \mathcal{H}_{d}/\{0\}}{\int_0^1(ra'^2+{d^2\over r}a^2)dr\over{1\over\ep^2}\int_0^1r(1-f_d^2({r\over\ep}))a^2dr}
\end{equation}
It is a classical result that these infimum are attained. Considering the rescaling  $(\tilde{a},\tilde{b})(r)=(a(\ep r),b(\ep r))$  and an extension by 0 outside $[0,1/\ep]$, we see that  $\ep\mapsto m_{\gamma_1,\gamma_2}(\ep)$ decreases when $\ep$ decreases. Then $\lim_{\ep\rightarrow0}m_{\gamma_1,\gamma_2}(\ep)$ exists.\\
Moreover, ${m}_{\gamma_1,\gamma_2}(\ep)$ is a simple eigenvalue and there exists an eigenvector $(a,b)$ verifying
$$a(r)\geq -b(r)\geq0\hbox{  for all $r>0$}.$$
Also, $m_0(\ep)$ is realized by some function $a(r)\geq0$.\\
In the previous works on the subject, the eigenvalue problem was  $\mathcal{L}_{n\mathcal{R}}\omega=-\lambda(\ep)\omega.$ We have
$$(\exists \lambda<0)\Leftrightarrow (\exists \mu<1).$$
By examining the proof of Theorems \ref{solutionscontinuesen0},
the possible behaviors at 0 of the solutions of the system
    \begin{equation}\label{eq:GLvptilde}
\left\{\begin{array}{rl}
a''+{a'\over r}-{\gamma_1^2\over r^2}a-f^2_da-f^2_db&=-\mu(\ep)(1-f^2_d)a\\
b''+{b'\over r}-{\gamma_2^2\over r^2}b-f^2_db-f^2_da&=-\mu(\ep)(1-f^2_d)b.
\end{array}
\right.
\end{equation}
and those of the solutions of (\ref{eq:GL}) are the same.
More precisely, if $\mu(\ep)$ is a bounded eigenvalue, it behaves as an
additional bounded parameter and we construct two solutions denoted by $(\alpha_1,\beta_1)$ and $(\alpha_3,\beta_3)$, depending on $\ep$ and verifying, for $r\in[0,R]$, 
$$\vert \alpha_1(r)\vert+\vert \beta_1(r)-r^{\gamma_2}\vert\leq Cr^{\gamma_2+2d+1}, \quad\vert \alpha'_1(r)\vert +\vert \beta'_1(r)-\gamma_2 r^{\gamma_2-1}\vert\leq Cr^{\gamma_2+2d}$$
 and
$$\vert \alpha_3(r)-r^{\gamma_1}\vert+\vert \beta_3(r)\vert\leq C r^{\gamma_1+2},\quad
 \vert \alpha'_3(r)-\gamma_1 r^{\gamma_1-1}\vert+\vert \beta'_3(r)\vert\leq Cr^{\gamma_1+1},$$
where $R$ and $C$ are independent of $\ep$, as in the proof of Theorem \ref{solutionscontinuesen0}.\\
We can suppose that 
 $\mu(\ep)\rightarrow\mu$, as $\ep\rightarrow0.$  
   Let 
  $\omega_{\ep}=(a_{\ep},b_{\ep})$
   be an eigenvector associated to $\mu(\ep)$. We define  
   $\tilde{\omega}_{\ep}(r)=\omega_{\ep}(\ep r),$ for $r\in[0,{1\over\ep}]$.
For some constants $A_{\ep}$ and $B_{\ep}$, 
 $(\tilde{a}_{\ep},\tilde{b}_{\ep})= A_{\ep}(\alpha_1,\beta_1)+B_{\ep}(\alpha_2,\beta_2).$
 We may suppose that $\max\{\vert A_{\ep}\vert,\vert B_{\ep}\vert\}=1$. Thus, $(\tilde{a}_{\ep},\tilde{a'}_{\ep},\tilde{b}_{\ep},\tilde{b'}_{\ep})(R)$ is bounded independently of $\ep$. Considering it as a Cauchy data, in the range $r\geq R$, we deduce that $(\tilde{a}_{\ep},\tilde{a'}_{\ep},\tilde{b}_{\ep},\tilde{b'}_{\ep})$ is bounded independently of $\ep$, in every interval $[R,\alpha]$, $\alpha>0$. Finally, we deduce the existence of some $\omega_0$ such that
 $$\tilde{\omega}_{\ep}\rightarrow \omega_0,\quad \hbox{as $\ep\rightarrow0$},$$ 
 uniformly on each compact subset of $[0,+\infty]$, where $\omega_0=(a_0,b_0)$ verifies
 \begin{equation}\label{eq:GLmu}
\left\{\begin{array}{rl}
a_0''+{a_0'\over r}-{\gamma_1^2\over r^2}a_0-f_d^2a_0-f_d^2b_0&=-\mu(1-f_d^2)a_0\\
b_0''+{b_0'\over r}-{\gamma_2^2\over r^2}b_0-f_d^2b_0-f_d^2a_0&=-\mu(1-f_d^2)b_0
\end{array}
\right.
\end{equation} 
Examining the proof of Theorem \ref{solutions continues en infty}, the possible behaviors at $+\infty$ of the solutions of (\ref{eq:GLmu}) are those given in Theorem \ref{solutions continues en infty}, when we suppose that ${\gamma^2_1+\gamma^2_2 \over2}-\mu d^2 >0$ and when we replace $n$ by $\sqrt{{\gamma^2_1+\gamma^2_2 \over2}-\mu d^2 }$.\\
Let us remark that the function $f_d$ and the eigenvalue problem used here are not exatly the same as in the previous works \cite{Miro}, \cite{TCL1} and \cite{TCL2} and \cite{SR}. However,
the proofs of the three following Theorems can be deduced from these works.
\begin{theorem} \label{lambda0}
For all $d\geq1$,\\
(i) there exist $C>0$ and $\ep_0>0$ such that, for all $\ep<\ep_0$,  ${m_{0}(\ep)-1\over\ep^2}\geq C$; $m_0(\ep)\rightarrow 1$ and there exists an associated eigenvector $a_{\ep}$ such that $\tilde{a}_{\ep}\rightarrow f_d$, uniformly in each $[0,R]$, $R>0$.\\
(ii) $m_{d-1,d+1}(\ep)>1$  and  ${m_{d-1,d+1}(\ep)-1\over\ep^2}\rightarrow 0$. \\
(iii) for $d>1$ and $n\geq 2d-1$, there exist $C>0$ and $\ep_0>0$ such that, for all $\ep<\ep_0$, ${m_{\vert d-n\vert,d+n}(\ep)-1\over\ep^2}\geq C$.\\
(iv) There exists an eigenvector $\omega_{\ep}$ associated to the eigenvalue $m_{d-1,d+1}(\ep)$  such that $\|(1-f^2_d)^{1\over2}(\tilde{\omega}_{\ep}-F_d)\|_{L^2(B(0,{1\over\ep}))}\rightarrow0$, as $\ep\rightarrow0$, where $F_d=(f_d'+{d\over r}f_d, f_d'-{d\over r}f_d)$ appears in Theorem \ref{essentiel}.\\
\end{theorem}
 The interested reader can find a direct proof of Theorem \ref{lambda0} in the appendix of \cite{HAL}.\\
 The following theorem is very important for our proof.
 \begin{theorem}\label{valeurproprenegative}
 Let $d\in\R$, $d>1$ be given. For all $n\in]1,d+1[$, there exists $C_n>0$ independent of $\ep$ such that
 $$m_{\vert d-n\vert,d+n}(\ep)\leq 1-C_n.$$
\end{theorem}
For the sake of completeness, we give a proof of this theorem in Part VI of the present paper, following the proof of  \cite{Miro}, given for $n=d=2$.\\
The following theorem connects the eigenvalue problem to the existence of the bounded solutions.
\begin{theorem}\label{lierlespb}
(i) Let  $d>0$ and $\gamma_2>\gamma_1\geq0$ be given. 
With the notation above, if $\mu(\ep)\rightarrow \mu$, if $\tilde{\omega}_{\ep}\rightarrow \omega_0$, if ${\gamma^2_2+\gamma^2_1\over2}-\mu d^2>0$ and if $\omega_0$ blows up at $+\infty$, then ${\mu(\ep)-1\over\ep^2}\geq C$, where $C$ is a given positive number, independent of $\ep$.\\
(ii) If there exists some bounded solution $(a,b)$ of (\ref{eq:1}), then there exists an eigenvalue $\mu(\ep)$ verifying
${\mu(\ep)-1}\rightarrow0$.
\end{theorem}
To make the paper as self contained as possible, we give the proof of Theorem \ref{lierlespb} (i) in Part VI, following the proof of \cite{TCL2}, given for $\mu=1$ and for the eigenvalue $\lambda(\ep)$. The interested reader can find a direct proof of Theorem \ref{lierlespb} (ii) in \cite{HAL}.\\
The following theorems are new.
\begin{theorem}\label{lapremiere}
 When ${\gamma^2_1+\gamma^2_2\over\ep^2}-d^2>0$, if there exists some bounded solution $\omega=(a,b)$ of (\ref{eq:1}), then
we have  $\lim_{\ep\rightarrow0}m_{\gamma_1,\gamma_2}(\ep)\geq1$.
 \end{theorem}
Combining Theorem \ref{lapremiere} and Theorem \ref{lierlespb} (ii), we get the following
 \begin{corollary}\label{solutionbornee}
 If there exists some bounded solution $\omega=(a,b)$ of (\ref{eq:1}), then
we have  $\lim_{\ep\rightarrow0}m_{\gamma_1,\gamma_2}(\ep)=1$ and if $\omega_{\ep}$ is some eigenvector associated to 
$m_{\gamma_1,\gamma_2}(\ep)$, then $\tilde{\omega}_{\ep}$ tends to $\omega$, uniformly in all $[0,R]$, $R>0$.
 \end{corollary}
 The following theorem can be deduced at once from Theorem \ref{lapremiere} and Theorem \ref{valeurproprenegative}. 
  \begin{theorem}\label{n<d+1}
Let $n$ and $d$ be real numbers and $\gamma_1=\vert n-d\vert$, $\gamma_2=n+d$. There is no bounded solution of (\ref{eq:1}), when $d\geq 1$ and $1<n<d+1$.
 \end{theorem}
Using Theorems \ref{essentiel}, \ref{lesextremes} and \ref{n<d+1}, we will prove the following theorem.
\begin{theorem}\label{ngeqd+1}
There is no bounded solution of (\ref{eq:1}), whenever $d\geq 1$ and $n\geq d+1$.
 \end{theorem}
Then Theorem \ref{pasdesolutionbornee} is proved. With Theorem \ref{lierlespb} (i), we get
\begin{theorem}\label{vptendant}
For $d\geq1$, $n>1$, $\gamma_1=\vert n-d\vert$ and $\gamma_2=n+d$,  there is no eigenvalue $\mu(\ep)$, with eigenvector in $\mathcal{H}_{\vert n-d\vert}$,  such that 
${\mu(\ep)-1\over\ep^2}\rightarrow0$, as $\ep\rightarrow0$.
\end{theorem}
The paper is organised as follows.
In Part II, we give a sketch of the proofs of Theorem \ref{fd}, Theorem \ref{solutionscontinuesen0} and Theorem \ref{solutions continues en infty}. Complete proofs of  Theorem \ref{solutionscontinuesen0} and Theorem \ref{solutions continues en infty} are altogether long,  technical and classical. The interested reader can consult Part II and Part III of \cite{HAL}, which is a long preliminary version of the present paper. In Part III, we prove Theorem \ref{lesextremes}. In Part IV,  we prove Theorem \ref{lapremiere}. In Part V, we prove Theorem \ref{ngeqd+1}. In Part VI, we give the proof of Theorem \ref{lierlespb} (i) and of Theorem \ref{valeurproprenegative}, which  is needed in the proof of Theorem \ref{n<d+1}. Theorem \ref{lierlespb} (i) is needed to prove Theorem \ref{lapremiere} and to deduce Theorem \ref{vptendant}.  \\
\section{Proof of Theorem \ref{fd}, proof of Theorem \ref{solutionscontinuesen0}, proof of Theorem \ref{solutions continues en infty}.}
\subsection{ Proof of Theorem \ref{fd}.}
 The existence of $f_d$, its expansion near 0 and $+\infty$ and its property of uniqueness are proved in \cite{HH}. However, these authors suppose that $d\in\N^{\star}$ and this is used only in the first step of their proof. Let us 
give an alternative proof for this first step, valid for all $d>0$. We have to prove that for all $a>0$ there exists some solution of (\ref{eq:f}) verifying $f\sim_0 ar^d$ and that $f$ is defined in an interval $[0,R]$, with $R>0$. We rewrite the equation (\ref{eq:f}) as
$$(r^{2d-1}(r^{-d}f)')'=-r^{2d-1}f(1-f^2).$$
For all $R>0$ and all $a>0$, $f$ solves (\ref{eq:f}) in $[0,R]$, and $f\sim_0 ar^d$ if and only if the map $g:r\mapsto r^{-d}f(r)$ is a fixed point in 
$\mathcal{C}([0,R])$ of the function $\Phi$ defined by
\begin{equation}\label{eq:Phi}
\Phi(g)(r)=a+\int_0^rt^{-2d+1}\int_0^t-s^{3d-1}g(s)(1-s^{2d}g^2(s))ds.
\end{equation}
Let us denote $\varphi(s,g)=-g(1-s^{2d}g^2)$. As in the proof of the Cauchy-Lipschitz Theorem, we remark first that for all $\alpha>0$ and all $\beta>0$ exist $M$ and $C$ such that
$$\left(s\in]0,\alpha],\quad\|g-a\|_{L^{\infty}([0,\alpha])}<\beta\right)\Rightarrow\left( \|\varphi(s,g)\|_{L^{\infty}([0,\alpha])}\leq M\right)$$
and
$$\left(s\in]0,\alpha],\quad\|g_1-a\|_{L^{\infty}([0,\alpha])}<\beta,\quad \|g_2-a\|_{L^{\infty}([0,\alpha])}<\beta\right)$$
$$\Rightarrow\left(\vert\varphi(s,g_1)-\varphi(s,g_2)\vert(s)\leq C\vert g_1-g_2\vert(s)\right).
$$
Moreover, 
$M$ and $C$ remain inchanged if $\alpha$ is replaced by a smallest positive number.\\
Now, we estimate, for $r\in[0,\alpha]$
$$\vert\Phi(g)-a\vert(r)\leq M{r^{d+2}\over 3d(d+2)}\hbox{  
and  }
\vert\Phi(g_1)-\Phi(g_2)\vert(r)\leq C{r^{d+2}\over 3d(d+2)}\|g_1-g_2\|_{L^{\infty}([0,\alpha])}.$$
Now, we choose some $R$ such that
$$0<R<\min\{1,\alpha,{3d(d+2)\beta\over M},{3d(d+2)\over C}\}$$
and we denote $\overline{B}(a,\beta)=\{g\in \mathcal{C}[0,R]);\|g-a\|_{L^{\infty}([0,R])}\leq \beta\}$,
in order $\Phi$ to be a contractant function from the closed subset $\overline{B}(a,R)$ of the Banach space $\mathcal{C}([0,R])$ into itself. Thus, by the Banach fixed Point Theorem, $\Phi$ has a unique fixed point $g$ in $\overline{B}(a,R)$.
Then $r\mapsto r^d g(r)$, defined in $[0,R]$, is the desired solution of (\ref{eq:f}).\\
The proof of \cite{HH} can be used to conclude to the existence of $A_d$. \\
Now, let us prove the continuity of $d\mapsto g^{(d)}$, where
 $g^{(d)}(r)=r^{-d}f_d(r)$. First, let us prove that the map $d\mapsto A_d$, defined in $\R^{+\star}$, increases.\\
As a first step, for $\delta\neq d$, we combine the equations of $f_d$ and $f_{\delta}$ to obtain, for every $(r_1,r_2)$, $0<r_1<r_2$, 
 $$[r(f'_df_{\delta}-f'_{\delta}f_d)]^{r_2}_{r_1}=\int^{r_2}_{r_1}f_df_{\delta}(f^2_d-f^2_{\delta})dt.$$
We derive two properties. The first one is that $f_d-f_{\delta}$ cannot keep the same sign in $[0,+\infty[$, otherwise, when $r_1=0$ and $r_2\rightarrow+\infty$, the lrs would be 0 and the rhs would be non zero. The second one is that $f_d$ and $f_{\delta}$ can be equal only for one value $r>0$. Indeed, if $r_1<r_2$ are such that $f_d(r_i)=f_{\delta}(r_i)$, for $i=1,2$, we get that
$ r_2f_d(r_2)(f_d-f_{\delta})'(r_2)-r_1f_d(r_1)(f_d-f_{\delta})'(r_1)$ has the same sign as $f^2_d-f^2_{\delta}$ in $[r_1,r_2]$, and this is a contradiction.
\\
Now, let $0<\delta<d$ be given. Near $+\infty$ we have the expansion 
$f_d(r)-f_{\delta}(r)={\delta^2-d^2\over 2r^2}+0({1\over r^4})$
and consequently, there exists $R>0$ such that $f_d(r)<f_{\delta}(r)$, for all $r\in [R,+\infty[$. But we have also $r^d<r^{\delta}$ for $0<r<1$. Since the sign of $f_d-f_{\delta}$ has to change once in $[0,+\infty[$, and in view of the expansions near 0, we deduce that $A_d>A_{\delta}$.\\
 Now, we denote $\lim_{d\rightarrow\delta,d>\delta}A_d=B$.
 But  
 $f_d$ is defined in $[0,+\infty[$. We have, for all $r>0$,
 $g^{(d)}=\Phi(g^{(d)})$, where $\Phi$ is defined in (\ref{eq:Phi}), but with $A_d$ instead of $a$. Using in addition $0\leq f_d(1-f^2_d)\leq1$, we get  that for all $\alpha>0$, there exists $\beta>0$ independent of $d$ in an interval containing $\delta$, such that $\vert g^{(d)}\vert(r)\leq \beta$ and $\vert(g^{(d)})'\vert(r)\leq\beta$, for all $r\in[0,\alpha]$. So, for all $r>0$, $g^{(d)}(r)$ has a limit, denoted by $g$, as $d\rightarrow\delta$, uniformly in every $[0,\alpha]$, $\alpha>0$ and we have $\Phi(g)(r)=g(r)$, for all $r>0$, where $\Phi$ is defined in (\ref{eq:Phi}), but with $B$ instead of $a$ and $\delta$ instead of $d$. Consequently, if we denote $f(r)=r^dg(r)$, then $f\sim_0 Br^{\delta}$, $f$ is a solution of (\ref{eq:f}) (with $\delta$ in place of $d$), $f$ is non decreasing in $[0,+\infty[$. In view of the uniqueness of such a solution of (\ref{eq:f}) (\cite{HH}), we deduce that $B=A_{\delta}$ and that $f=f_{\delta}$. The same result remains true when $d\rightarrow\delta$, $d<\delta$. We have proved that $d\mapsto g^{(d)}$ is continuous from $[1,+\infty[$ into $L^{\infty}([0,\alpha])$, for all $\alpha>0$. 
\subsection{Sketch of the proof of Theorem \ref{solutionscontinuesen0}. }
The pattern of proof is the same for the four solutions. Let us give an idea of the proof.

1.  We construct some solution $(a_1,b_1)$ such that for all compact subset $\mathcal{K}$ of $\mathcal{D}$, there exist some $R>0$, depending only on $\mathcal{K}$ and some $C>0$, also depending only on $\mathcal{K}$, such that for all $r\in]0,R]$ and all $(d,\gamma_1,\gamma_2)\in\mathcal{K}$, we have
$$\vert a_1(r)\vert+\vert b_1(r)-r^{\gamma_2}\vert\leq Cr^{\gamma_2+2d+1}\quad\hbox{
and}\quad
\vert a'_1(r)\vert +\vert b'_1(r)-\gamma_2 r^{\gamma_2-1}\vert\leq Cr^{\gamma_2+2d}.$$
and such that, for all $r\in]0,R]$,  $(d,\gamma_1,\gamma_2)\mapsto (a_1(r),a'_1(r),b_1(r),b'_1(r))$ is continuous on $\mathcal{K}$, and differentiable wrt $\gamma_1$ and wrt $\gamma_2$. First, the construction is done for $r\in]0,R]$. Then the definition of this solution in $[0,+\infty[$ and the continuity wrt $(d,\gamma_1,\gamma_2)\in\mathcal{K}$, for all $r>0$, follows from the Cauchy-Lipschitz Theorem. Let us remark the importance for the constants $C$ and $R$ to be independent of the parameters. \\
We use a constructive method, similar to the proof of the Banach fixed point Theorem.\\
We define a fixed point problem of the form
$(a,b)=\Phi(a,b),$
that is 
\begin{equation}\label{eq:omega1}
\left\{\begin{array}{rl}
a&=r^{\gamma_1}+r^{\gamma_1}\int_0^rt^{-2\gamma_1-1}\int_{0}^ts^{\gamma_1+1}(f_d^2b-(1-2f_d^2)a)dsdt\\
b&=r^{\gamma_2}\int_0^rt^{-2\gamma_2-1}\int_0^ts^{\gamma_2+1}(f_d^2a-(1-2f_d^2)b)dsdt.
\end{array}
\right.
\end{equation}
whose solutions verify the differential system that we have to solve. 

2. We construct some solution $(a_3,b_3)$, such that, for any compact subset  $\mathcal{K}\in\mathcal{D}$, exist some real numbers $R$ and $C$ verifying, for all $0<r< R$,
$$\vert a_3(r)-r^{\gamma_1}\vert\leq C r^{\gamma_1+2},\vert b_3(r)\vert\leq Cr^{\gamma_1+2d+2},
\vert a'_3(r)-\gamma_1 r^{\gamma_1-1}\vert\leq Cr^{\gamma_1+1},\vert b'_3(r)\vert\leq C r^{\gamma_1+2d+1}.$$
 For this purpose, we consider the fixed point problem  
\begin{equation}\label{eq:omega3}
\left\{\begin{array}{rl}
a&=r^{\gamma_1}\int_0^rt^{-2\gamma_1-1}\int_{0}^ts^{\gamma_1+1}(f_d^2b-(1-2f_d^2)a)dsdt\\
b&=r^{\gamma_2}+r^{\gamma_2}\int_0^rt^{-2\gamma_2-1}\int_0^ts^{\gamma_2+1}(f_d^2a-(1-2f_d^2)b)dsdt.
\end{array}
\right.
\end{equation}

3. For the construction of $(a_2,b_2)$, in the case  when $(d,\gamma_1,\gamma_2)\in\mathcal{D}_1$, we consider the fixed point problem
\begin{equation}\label{eq:gamma2}
 \left\{\begin{array}{rl}
a&= r^{-\gamma_1}\int_0^rt^{2\gamma_1-1}\int_{1}^ts^{-\gamma_1+1}(f_d^2b -(1-2f_d^2)a)dsdt\\
b&=r^{-\gamma_2}+r^{-\gamma_2}\int_0^rt^{2\gamma_2-1}\int_1^ts^{-\gamma_2+1}(f_d^2a-(1-2f_d^2)b)dsdt.
\end{array}
\right.
 \end{equation}
 while, when $(d,\gamma_1,\gamma_2)\in\mathcal{D}_2$, we consider the fixed point problem
\begin{equation}\label{eq:gamma2bis}
 \left\{\begin{array}{rl}
a&= r^{\gamma_1}\int_0^rt^{-2\gamma_1-1}\int_{0}^ts^{\gamma_1+1}(f_d^2b -(1-2f_d^2)a)dsdt\\
b&=r^{-\gamma_2}+r^{-\gamma_2}\int_0^rt^{2\gamma_2-1}\int_1^ts^{-\gamma_2+1}(f_d^2a-(1-2f_d^2)b)dsdt.
\end{array}
\right.
 \end{equation}

4. In order to construct a solution $(a_4,b_4)$, when $(d,\gamma_1,\gamma_2)\in\mathcal{D}_1$, we solve 
 the following fixed point problem
 \begin{equation}\label{eq:omega4}
 \left\{\begin{array}{rl}
a&=r^{-\gamma_1}+ r^{-\gamma_1}\int_0^rt^{2\gamma_1-1}\int_{1}^ts^{-\gamma_1+1}(f_d^2b -(1-2f_d^2)a)dsdt\\
b&=r^{-\gamma_2}\int_0^rt^{2\gamma_2-1}\int_1^ts^{-\gamma_2+1}(f_d^2a-(1-2f_d^2)b)dsdt.
\end{array}
\right.
 \end{equation}
and, when $(d,\gamma_1,\gamma_2)\in\mathcal{D}_2$, we solve 
 the following fixed point problem
\begin{equation}\label{eq:gamma4bis}
\left\{\begin{array}{c}
a=\tau(r)+\tau(r)\int_0^r{1\over t}\tau^{-2}(t)\int_0^ts\tau(s)(f_d^2 b-(1-2f_d^2)a)dsdt\\
b=r^{-\gamma_2}\int_0^rt^{2\gamma_2-1}\int_0^t s^{-\gamma_2+1}(f_d^2 a-(1-2f_d^2)b)dsdt
\end{array}
\right.
\end{equation}
\subsection{Sketch of the proof of Theorem \ref{solutions continues en infty}.}
We use the system (\ref{eq:2}) and we construct a base of four solutions, $(x_j,y_j)$, $j=1,\ldots, 4$, characterized by their behaviors at $+\infty$. The solutions $(u_j,v_j)$ announced in Theorem \ref{solutions continues en infty} are obtained by  $u_j={x_j+y_j\over2}$ and $v_j={x_j-y_j\over2}$.\\
We denote
$$J_+={e^{\sqrt 2 r}\over\sqrt r},\quad J_-={e^{-\sqrt 2 r}\over\sqrt r},\quad\gamma^2={\gamma^2_1+\gamma^2_2\over2},\quad n=\sqrt{\gamma^2-d^2},\quad\xi^2={\gamma^2_2-\gamma^2_1\over2}.$$
We can replace the first equation of 
(\ref{eq:2}) by
$$(e^{2\sqrt 2 r}({x}e^{-\sqrt 2 r})')'=e^{\sqrt 2 r}q(r){x}-{\xi^2\over r^2} y\hbox{   
  or  }
(e^{-2\sqrt 2 r}({x}e^{\sqrt 2 r})')'=e^{-\sqrt 2 r}q(r){x}-{\xi^2\over r^{2}} y,$$
where
$$q(r)={-\gamma^2-3d^2\over r^2}+3(1-f_d^2+{d^2\over r^2}).$$
The second equation of the system (\ref{eq:2}) can be written as
$$(r^{2n+1}(r^{-n}y)')'=r^{n+1}({\xi^2\over r^{2}}x-(1-f_d^2-{d^2\over r^2})y)$$
or
$$(r^{-2n+1}(r^{n}y)')'=r^{-n+1}({\xi^2\over r^{2}}x-(1-f_d^2-{d^2\over r^2})y).$$
Finally, the system (\ref{eq:2})  can be written as
\begin{equation}\label{eq:systenlinfini}
\left\{\begin{array}{ll}
(e^{\pm 2\sqrt 2 r}(r^{1\over2}e^{\mp\sqrt 2 r}x)')'=r^{1\over2}e^{\pm\sqrt 2 r}q(r)x-{\xi^2\over r^2} y\\
(r^{\pm 2n+1}(r^{\mp n}y)')'=r^{\pm n+1}({\xi^2\over r^{2}}x-(1-f_d^2-{d^2\over r^2})y)
\end{array}
\right.
\end{equation}
In order to construct four solutions of (\ref{eq:systenlinfini}), we give $R_0>0$ and we define fixed points problems of the form
$(x,y)=\Phi(x,y),$
for $(x,y)$ defined in $[R_0,+\infty[$, and whose solutions are solutions of (\ref{eq:systenlinfini}). The function $\Phi$ will depend on $R_0$, except for one solution denoted by $(x_2,y_2)$ (vanishing exponentially at $+\infty$). The present construction does not allow us 
to construct the solutions $(x_j,y_j)$, $j\neq2$ without taking into account a given compact subset 
\begin{equation}\label{eq:defK}
\mathcal{K}\subset \{(d,\gamma_1,\gamma_2);0\leq\gamma_1<\gamma_2;\xi^2-d^2>0\}.
\end{equation}
 Indeed, $R_0$ depends on $\mathcal{K}$. Let us list the different fixed point problems we need.
 
1. The exponential blowing up behavior at $+\infty$ : the solution $(x_1,y_1)$. For $R_0>0$
$$\left\{\begin{array}{ll}
x&=J_++J_+\int_{+\infty}^r(J_+)^{-2}{1\over t}\int_{R_0}^tsJ_+({\xi^2\over s^2}y-3(1-f^2_d-{d^2\over s^2})x)dsdt\\
y&=r^n\int_{R_0}^rt^{-2n-1}\int_{R_0}^ts^{n+1}({\xi^2\over s^2}x-(1-f^2_d-{d^2\over s^2})y)dsdt.
\end{array}
\right.
$$

2. The intermediate blowing up behavior at $+\infty$ : the solution $(x_3,y_3)$. For $R_0>0$ 
$$\left\{\begin{array}{ll}
x&=J_+\int_{+\infty}^r(J_+)^{-2}{1\over t}\int_{R_0}^tsJ_+({\xi^2\over s^2}y-3(1-f^2_d-{d^2\over s^2})x)dsdt\\
y&=r^n+r^n\int_{+\infty}^rt^{-2n-1}\int_{R_0}^ts^{n+1}({\xi^2\over s^2}x-(1-f^2_d-{d^2\over s^2})y)dsdt
\end{array}
\right.
$$

3. The least behavior at $+\infty$ : the solution $(x_2,y_2)$.
We consider 
$$\left\{\begin{array}{ll}
x&=J_-+J_-\int_{+\infty}^r(J_-)^{-2}{1\over t}\int_{+\infty}^tsJ_-({\xi^2\over s^2}y-3(1-f^2_d-{d^2\over s^2})x)dsdt\\
y&=r^{-n}\int_{+\infty}^rt^{2n-1}\int_{+\infty}^ts^{-n+1}({\xi^2\over s^2}x-3(1-f^2_d-{d^2\over s^2})y)dsdt
\end{array}
\right.
$$

4. The intermediate vanishing behavior at $+\infty$ : the solution $(x_4,y_4)$.
For $R_0>0$ 
$$\left\{\begin{array}{ll}
x&=J_-\int_{R_0}^r(J_-)^{-2}{1\over t}\int_{+\infty}^tsJ_-({\xi^2\over s^2}y-3(1-f^2_d-{d^2\over s^2})x)dsdt\\
y&=r^{-n}+r^{-n}\int_{+\infty}^rt^{2n-1}\int_{+\infty}^ts^{-n+1}({\xi^2\over s^2}x-3(1-f^2_d-{d^2\over s^2})y)dsdt
\end{array}
\right.
$$
We need the following estimate, which is not difficult to prove, by an integration by part.
Let $\alpha\in\R$ and $\beta>0$ be given. Then
\begin{equation}\label{eq:J-}
\int_t^{+\infty}s^{\alpha}e^{-\beta s}ds\leq{2\over\beta}t^{\alpha}e^{-\beta t}\quad\hbox{  for all  }t\geq{2\alpha\over\beta}
\end{equation}
and 
\begin{equation}\label{eq:J+}
\int_R^{t}s^{\alpha}e^{\beta s}ds\leq{2\over\beta}t^{\alpha}e^{\beta t}\quad\hbox{  for all  }t\geq R\geq{-2\alpha\over\beta}\quad
\end{equation}
 \section{The smallest behavior at zero is connected with the greatest behavior at infinity. }
 {\bf Proof of Theorem \ref{lesextremes}.}\\
Let $(d,\gamma_1,\gamma_2)\in\mathcal{D}$.
 Let us prove first that $(a_1,b_1)$ blows up exponentially at $+\infty$.\\
Let us define $x=a_1+b_1$ and $y=a_1-b_1$. We have $x(r)\sim_0 r^{\gamma_2}$ and $y(r)\sim_0 -r^{\gamma_2}$. Thus, we have $x(r)>0$ and $y(r)<0$ near $r=0$.
 Let us suppose that $x(r)>0$ and $y(r)<0$ in $]0,R[$. Combining the first equation of the system (\ref{eq:2}) and the equation (\ref{eq:f}), we get, for all $r\geq0$
$$[rx'f_d-rf_d'x]_0^{r}+\int_0^{r}{-\gamma^2+d^2\over s}xf_dds+\mu^2\int_0^{r}{y\over s}f_d ds-2\int_0^{r} sf_d^3xds=0.$$
For $0<r\leq R$, we deduce that
\begin{equation}\label{eq:xincreases}
rf_d^2({x\over {f_d}})'({r})\geq 2\int_0^{r} sf_d^3xds.
\end{equation}
This proves that ${x\over {f_d}}$ increases in $]0,R]$ and therefore $x(R)>0$.\\
Moreover, combining the second equation of the system (\ref{eq:2}) and (\ref{eq:f}), we get
$$[ry'f_d-rf_d'y]_0^{r}+\int_0^{r}{-\gamma^2+d^2\over s}yf_dds+\xi^2\int_0^{r}{x\over s}f_d ds=0.$$
For $0<r\leq R$, we deduce that
\begin{equation}\label{eq:ydecreases}
rf_d^2({-y\over {f_d}})'({r})\geq \int_0^{r}{-\gamma^2+d^2\over s}yf_dds.
\end{equation}
This proves that ${-y\over{f_d}}$ increases in $]0,R]$ and therefore $-y(R)>0$. 
Finally, we have proved that $x(r)>0$ and $y(r)<0$ for all $r>0$. Now (\ref{eq:xincreases}) and (\ref{eq:ydecreases}) are valid for all $r>0$ and we know that $f_d\sim_{+\infty} 1$. Thus, the behavior of $x$ at $+\infty$ cannot be a polynomial increasing behavior. We return to Theorem \ref{solutions continues en infty} that gives all the possible behaviors at $+\infty$ and we deduce that $x$ and $y$ have an exponentially increasing behavior at $+\infty$. So $a$ and $b$ have an exponentially increasing behavior at $+\infty$, too. 
\\
Let us prove now that $(u_2,v_2)\sim_0D (o(r^{\gamma_1}),r^{-\gamma_2})$, for some $D\neq0$. Multiplying (\ref{eq:1}) and integrating by parts, we get easily, for all $r_1>0$ and $r>0$
 $$[s(a'_1u_2-u'_2a_1+v_2b'_1-v'_2b_1)(s)]_{r}^{r_1}=0.$$
Using $(a_1,b_1)\sim_{+\infty} C({e^{\sqrt2 r}\over\sqrt r},{e^{\sqrt2 r}\over\sqrt r})$, for some $C\neq0$, and $(u_2,v_2)\sim_{+\infty}({e^{-\sqrt2 r}\over\sqrt r},{e^{-\sqrt2 r}\over\sqrt r})$, we get 
$$\lim_{r\rightarrow+\infty}r(a'_1u_2-u'_2a_1+v_2b'_1-v'_2b_1)(r)=4C\sqrt 2.$$
 Consequently
$$\lim_{r\rightarrow 0}r(a'_1u_2-u'_2a_1+v_2b'_1-v'_2b_1)(r)=4C\sqrt 2.$$
We know that $(a_1,b_1)\sim_0 (o(r^{\gamma_2}),r^{\gamma_2})$. According to Theorem \ref{solutionscontinuesen0}, that gives all the possible behaviors at 0, we conclude that the only fitting behavior at 0 for $(u_2,v_2)$ is $(u_2,v_2)\sim_0 D (o(r^{\gamma_1}),r^{-\gamma_2})$, for $D={2C\sqrt2\over \gamma_2}$.\\
This ended the proof of Theorem \ref{lesextremes}.
\section{The proof of Theorem \ref{lapremiere} and of Corollary \ref{solutionbornee}.}
 Let $d>1$. We can rewrite the system (\ref{eq:1}) as 
 \begin{equation}\label{eq:formematrice}
 X'=MX \hbox{  with  } X=(a,ra',b,rb')^{\hbox{t}}
  \end{equation}
with
 $$M=\left(\begin{array}{cccc}
 0&{1\over r}&0&0\\
 -r(1-2f_d^2)+{\gamma_1^2\over r}&0 &rf_d^2&0\\
 0&0&0&{1\over r}\\
 rf_d^2&0& -r(1-2f_d^2)+{\gamma_2^2\over r}&0
 \end{array}
 \right).
$$ 
\begin{lemma} \label{systemefondamental}
Let us suppose that there exists a bounded solution of (\ref{eq:1}) and let us chose a base of solutions,
 $X_1$, $X_2$, $X_3$, $X_4$,  for (\ref{eq:formematrice}), whose third vector is a bounded solution. Let us name  $R(s)$ the resolvant matrix, whose columns are the vectors $X_i$, $i=1,\ldots,4$. Let us name $\mathcal{C}_2$ and $\mathcal{C}_4$ the second and the fourth column of $R^{-1}(s)$. We have

 $$\hbox{  at   $0$ and  when  $(d,\gamma_1,\gamma_2)\in\mathcal{D}_1$ and $\gamma_1+\gamma_2-2d-2<0$    }$$
 $$\begin{array}{rl}\mathcal{C}_2=\left(\begin{array}{c}
 O(s^{\gamma_1})\\
O(s^{\gamma_1+2\gamma_2})\\
 O(s^{-\gamma_1})\\
 O(s^{\gamma_1})
 \end{array}
 \right)
\hbox{  and  } &\mathcal{C}_4=\left(\begin{array}{c} O(s^{-\gamma_2})\\
O(s^{\gamma_2})\\
O(s^{\gamma_2})\\
O(s^{2\gamma_1+\gamma_2})
\end{array}
\right)
\end{array}
 $$
 and
 $$\hbox{  at   $0$ and  when  $(d,\gamma_1,\gamma_2)\in\mathcal{D}_1$ and $\gamma_1+\gamma_2-2d-2>0$    }$$
 $$
 \begin{array}{rl}\mathcal{C}_2=\left(\begin{array}{c}
 O(s^{-\gamma_2+2d+2})\\
O(s^{\gamma_2+2d+2})\\
 O(s^{-\gamma_1})\\
 O(s^{\gamma_1})
 \end{array}
 \right)
\hbox{  and  } &\mathcal{C}_4=\left(\begin{array}{c} O(s^{-\gamma_2})\\
O(s^{\gamma_2})\\
O(s^{-\gamma_1+2d+2})\\
O(s^{\gamma_1+2d+2})
\end{array}
\right)
\end{array}
 $$
 and
 $$\hbox{  at   $0$ and  when  $(d,\gamma_1,\gamma_2)\in\mathcal{D}_2$    }$$
 $$\begin{array}{rl}\mathcal{C}_2=\left(\begin{array}{c}
 O(\tau(s)s^{-\gamma_2+\gamma_1+2d+2})\\
O(\tau(s) s^{\gamma_1+\gamma_2+2d+2})\\
 O(\tau(s))\\
 O(s^{\gamma_1})
 \end{array}
 \right)
\hbox{  and  } &\mathcal{C}_4=\left(\begin{array}{c} O(s^{\gamma_1-\gamma_2}\tau(s))\\
O(s^{\gamma_1+\gamma_2}\tau(s))\\
O(s^{2d+2}\tau(s))\\
O(s^{\gamma_1+2d+2})
\end{array}
\right)
\end{array}
 $$ 
 and in any case, at $+\infty$
 $$\begin{array}{rl}\mathcal{C}_2\sim_{+\infty}{1\over -16n\sqrt2}\left(\begin{array}{c}
 4nJ_-\\
4n J_+\\
-4\sqrt 2 s^{n}\\
-4\sqrt 2 s^{-n}
 \end{array}
 \right)
  \hbox{    and   }\quad \mathcal{C}_4\sim_{+\infty}{1\over -16n\sqrt2}\left(\begin{array}{c}
 4nJ_-\\
4n J_+\\
4\sqrt 2 s^{n}\\
4\sqrt 2 s^{-n}
 \end{array}
 \right)
 \end{array}
 $$
 where $-16n\sqrt2$ is the determinant of $R(s)$.
 \end{lemma}
  {\bf Proof } 
 $R(s)$ is chosen as follows
 $$R(s)\sim_{+\infty}\left(\begin{array}{cccc}
 J_+&J_-&s^{-n}&s^n\\
 s(J_+)'&s(J_-)'&-ns^{-n}&ns^n\\
  J_+&J_-&-s^{-n}&-s^n\\
 s(J_+)'&s(J_-)'&ns^{-n}&-ns^n\\
 \end{array}
 \right) 
 $$
 where, as usual, the notation $J_+$ stands for ${e^{\sqrt2 s}\over\sqrt s}$ and the notation $J_-$ stands for  ${e^{-\sqrt2 s}\over\sqrt s}$.\\
To give the behaviors at 0, we return to Theorem \ref{solutionscontinuesen0}.  We have, for some $c_i\neq0$, $i=1,\ldots,4$
$$\hbox{If $(d,\gamma_1,\gamma_2)\in\mathcal{D}_1$,}\quad R(s)\sim_{0}\left(\begin{array}{cccc}
 O(s^{\gamma_2+2d+2})&O(s^{\tilde{\gamma}_1})&c_3s^{\gamma_1}&c_4s^{-\gamma_1}\\
 O(s^{\gamma_2+2d+2})&O(s^{\tilde{\gamma}_1})&c_3\gamma_1s^{\gamma_1}&-c_4\gamma_1s^{-\gamma_1}\\
  c_1s^{\gamma_2}&  c_2s^{-\gamma_2}&O(s^{\gamma_1+2d+2})&O(s^{\tilde{\gamma}_2})\\
 c_1\gamma_2s^{\gamma_2}&-c_2\gamma_2  s^{-\gamma_2}&O(s^{\gamma_1+2d+2})&O(s^{\tilde{\gamma}_2})\\
 \end{array}
 \right) 
 $$
where we use the notation $$\tilde{\gamma}_1=\min\{\gamma_1,-\gamma_2+2d+2\}
\hbox{ and  }\tilde{\gamma}_2=\min\{\gamma_2,-\gamma_1+2d+2\}\quad\hbox{if }\gamma_1+\gamma_2-2d-2\neq0$$
(if $ \gamma_1+\gamma_2-2d-2=0$, we have to replace $O(s^{\tilde{\gamma}_1})$ by $O(s^{\gamma_1}\log s)$ and 
$O(s^{\tilde{\gamma}_2})$ by $O(s^{\gamma_2}\log s)$)

and 

 $$\hbox{If $(d,\gamma_1,\gamma_2)\in\mathcal{D}_2$,}\quad R(s)\sim_{0}\left(\begin{array}{cccc}
 O(s^{\gamma_2+2d+2})&O(s^{-\gamma_2+2d+2})&c_3s^{\gamma_1}&c_4\tau(s)\\
 O(s^{\gamma_2+2d+2})&O(s^{-\gamma_2+2d+2})&c_3\gamma_1s^{\gamma_1}&-c_4s\tau'(s)\\
  c_1s^{\gamma_2}&  c_2s^{-\gamma_2}&O(s^{\gamma_1+2d+2})&O(\tau(s)s^{2d+2})\\
 c_1\gamma_2s^{\gamma_2}&-c_2\gamma_2  s^{-\gamma_2}&O(s^{\gamma_1+2d+2})&O(\tau(s)s^{2d+2})\\
 \end{array}
 \right) 
 $$
where 

$$\tau(s)=\left\{\begin{array}{rl}{s^{-\gamma_1}-s^{\gamma_1}\over 2\gamma_1}\hbox{  if  }\gamma_1\neq0\\
-\log s\hbox{  if   }\gamma_1=0
\end{array}
\right.
$$

 The determinant $W$ of $R(s)$ is independent of $s$, due to the fact that the matrix $M$ of the differential system
has a null trace. Moreover, $J_+J_-={1\over s}$. Using the behavior at $+\infty$ of $R(s)$, given above, we deduce that $W$ is the principal term, as $s\rightarrow+\infty$ of 
$${1\over s}\left\vert\begin{array}{cccc}
1&1&1&1\\
s\sqrt2&-s\sqrt2&-n&n\\
1&1&-1&-1\\
s\sqrt2&-s\sqrt2&n&-n
\end{array}
\right\vert
$$
that is 
$$W=-16n\sqrt2.$$

A direct calculation of the suitable determinants gives the estimate of $\mathcal{C}_2$ and $\mathcal{C}_4$.\\

{\bf The proof of Theorem \ref{lapremiere} completed.}\\

Let $m=\lim_{\ep\rightarrow0}m_{\gamma_1,\gamma_2}(\ep)$.
We can define $\omega_{\ep}\in\mathcal{H}_{\gamma_1}$ an eigenvector associated to $m_{\gamma_1,\gamma_2}(\ep)$ and $\omega_0=(a_0,b_0)$ such that $\tilde{\omega}_{\ep}\rightarrow\omega_0$ on each compact subset of $[0,+\infty[$. 
In what follows, let us suppose that $m<1$. Then ${\gamma_1^2+\gamma_2^2\over2}-m d^2>0$.\\
 Since $a_0\geq-b_0\geq0$, the possible behaviors at $+\infty$ for $(a_0,b_0)$ are
 $(r^{-n_0},-r^{-n_0})\hbox{  and }(r^{n_0},-r^{n_0})$
 where 
 \begin{equation}\label{eq:n0}
 n_0=\sqrt{{\gamma_1^2+\gamma_2^2\over2}-md^2}.
 \end{equation}
 Since $m<1$, we have by Theorem \ref{lierlespb} (i), that $\omega_0$ has a bounded behavior at $+\infty$ and consequently
 $$(a_0,b_0)\sim_{+\infty} (r^{-n_0},-r^{-n_0})\quad\hbox{
 and}\quad a_0+b_0=O(r^{-n_0-2})\quad\hbox{at $+\infty$}.$$
At $0$, in view of $a_0\geq -b_0\geq0$, the only possible behavior is 
$$(a_0,b_0)\sim_0(cr^{\gamma_1}, O(r^{\gamma_1+2d+2})),\quad\hbox{
for some $c>0$}.$$
Let us prove that the hypothesis $m<1$ leads to a contradiction.\\
Since $n=\sqrt{{\gamma_1^2+\gamma_2^2\over2}-d^2}$, we have, by (4.34) 
$$(m<1)\Leftrightarrow (n_0>n).$$
 Let us denote 
$X_0=(
 a_0,
 ra'_0,
 b_0,
 rb'_0
 )^t
 $, the vector corresponding to $\omega_0$.
 We have 
 $$X'_0=MX_0-(m-1)(1-f^2_d)(0,ra_0,0,rb_0)^{\hbox{t}}.$$
 let us define $X_1$, $X_2$, $X_3$ and $X_4$ as in Lemma \ref{systemefondamental}. 
We are going to prove that there exist some constants $C_i$ such that
$$
 X_0=\sum_{i=1}^4 C_iX_i-(m-1)\sum_{i=1}^4\hat{ X}_i,$$
 with
  \begin{equation}\label{eq:X0-1}
 \hat{X}_i\hbox{ bounded at $0$, $i=1,2,3,4$}
  \end{equation}
and
  \begin{equation}\label{eq:X0-2} 
  \hbox{at $+\infty$}\left\{\begin{array}{c}
\hat{X}_1=X_1 O(r^{-n_0-3}J_-)  \quad;\quad\hat{X}_2=X_2O(r^{-n_0-3}J_+)\\
\hat{X}_3=X_3 O(1)\quad;\quad \hat{X}_4=X_4O(1).
 \end{array}
  \right.
 \end{equation}
In order to prove (\ref{eq:X0-1}) and (\ref{eq:X0-2}), we write
\begin{equation}\label{eq:sum}
 X_0=\sum_{i=1}^4A_i(r)X_i
 \end{equation}
with
 \begin{equation}\label{eq:Ai}
i=1,\ldots,4,\quad A_i(r)=A_i-(m-1)\int_1^r[R^{-1}(s)s(1-f_d^2)\left(\begin{array}{c}
0\\
a_0 \\
0\\
b_0
\end{array}
 \right)ds]_i
 \end{equation}
where  the notation $[\hbox{   }]_i$ means the $i^{\hbox{th}}$ line of the vector, and where $A_i$ is a constant.\\
  Let us examine the behavior of each term $A_i(r) X_i$ at $+\infty$ and at 0, using Lemma \ref{systemefondamental}.\\
For the first term, we use the first terms of $\mathcal{C}_2$ and $\mathcal{C}_4$, given in Lemma \ref{systemefondamental}, to obtain
$$[R^{-1}(s)s(1-f_d^2)\left(\begin{array}{c}
0\\
a_0 \\
0\\
b_0
\end{array}
 \right)]_1\sim_{+\infty}O({1\over s}J_-(a_0+b_0))
 $$
 $$\hbox{ and }\sim_0 \left\{\begin{array}{c}
 s(O(s^{\gamma_1}a_0+O(s^{-\gamma_2}b_0))\hbox{  if   }(d,\gamma_1,\gamma_2)\in\mathcal{D}_1,\gamma_1+\gamma_2-2d-2<0\\
 s(O(s^{-\gamma_2+2d+2}a_0+O(s^{-\gamma_2}b_0))\hbox{  if   }(d,\gamma_1,\gamma_2)\in\mathcal{D}_1,\gamma_1+\gamma_2-2d-2>0\\
 s(O(\tau(s)s^{\gamma_1-\gamma_2+2d+2}a_0+O(\tau(s)s^{\gamma_1-\gamma_2}b_0))\hbox{  if   }(d,\gamma_1,\gamma_2)\in\mathcal{D}_2.
 \end{array}
 \right.
 $$

Let us define
 $$B_1=-(m-1)\int_1^{+\infty}[R^{-1}(s)s(1-f_d^2)\left(\begin{array}{c}
0\\
a_0 \\
0\\
b_0
\end{array}
 \right)]_1ds\quad \hbox{ and  }\hat{X}_1=X_1\int_{+\infty}^r[R^{-1}(s)s(1-f_d^2)\left(\begin{array}{c}
0\\
a_0 \\
0\\
b_0
\end{array}
 \right)]_1ds$$
We can write
$$A_1(r)X_1=(A_1+B_1)X_1-(m-1)\hat{X}_1
$$
 We see that  $\hat{X}_1=X_1O(1)$ at 0.  Using (\ref{eq:J-}), we get
$\hat{X}_1=X_1O(r^{-n_0-3}J_-)$ at $+\infty$.

For the second term, we obtain
$$[R^{-1}(s)s(1-f_d^2)\left(\begin{array}{c}
0\\
a_0 \\
0\\
b_0
\end{array}
 \right)]_2\sim_{+\infty}O({1\over s}J_+(a_0+b_0))
 $$
$$
 \hbox{ and }\sim_0\left\{\begin{array}{c}
 s(O(s^{\gamma_1+2\gamma_2}a_0)+O(s^{\gamma_2}b_0))\hbox{  if  }(d,\gamma_1,\gamma_2)\in\mathcal{D}_1\hbox{  and  }\gamma_1+\gamma_2-2d-2<0\\
  s(O(s^{\gamma_2+2d +2}a_0)+O(s^{\gamma_2}b_0))\hbox{  if  }(d,\gamma_1,\gamma_2)\in\mathcal{D}_1\hbox{  and  }\gamma_1+\gamma_2-2d-2>0\\
 s\tau(s)(O(s^{\gamma_1+\gamma_2+2d+2})a_0+O(s^{\gamma_1+\gamma_2})b_0)\hbox{  if  }(d,\gamma_1,\gamma_2)\in\mathcal{D}_2
 \end{array}
 \right.
 $$

Denoting  $$B_2=-(m-1)\int_1^0[R^{-1}(s)s(1-f_d^2)\left(\begin{array}{c}
0\\
a_0 \\
0\\
b_0
\end{array}
 \right)]_2ds\quad\hbox{and }\hat{X}_2=X_2\int_{0}^r[R^{-1}(s)s(1-f_d^2)\left(\begin{array}{c}
0\\
a_0 \\
0\\
b_0
\end{array}
 \right)]_2ds$$
 we get
$$A_2(r)X_2=(A_2+B_2)X_2-(m-1)\hat{X_2}$$
with, by (\ref{eq:J+})
$$
\hat{X}_2=X_2O(r^{-n_0-3}J_+)\hbox{ at $+\infty$ }.$$
Moreover, $\hat{X}_2$ is bounded at $0$.\\

 For the third term, we obtain
 \begin{equation}\label{eq:3emetermeenlinfini}
[R^{-1}(s)s(1-f_d^2)\left(\begin{array}{c}
0\\
a_0 \\
0\\
b_0
\end{array}
 \right)]_3\sim_{+\infty}{-1\over 16n\sqrt2}{4\sqrt2 d^2\over s}s^{n}(-a_0+b_0).
 \end{equation}
 Since $-a_0+b_0\sim_{+\infty}-2r^{-n_0}$, then this term is integrable at $+\infty$.
 At $0$, it is
$$
\sim_0\left\{\begin{array}{c}
 s(O(s^{-\gamma_1}a_0)+O(s^{\gamma_2}b_0))\hbox{  if  }(d,\gamma_1,\gamma_2)\in\mathcal{D}_1\hbox{  and  }\gamma_1+\gamma_2-2d-2<0\\
  s(O(s^{-\gamma_1}a_0)+O(s^{-\gamma_1+2d+2}b_0))\hbox{  if  }(d,\gamma_1,\gamma_2)\in\mathcal{D}_1\hbox{  and  }\gamma_1+\gamma_2-2d-2>0\\
 s(O(\tau(s)a_0)+O(\tau(s)s^{2d+2}b_0))\hbox{  if  }(d,\gamma_1,\gamma_2)\in\mathcal{D}_2\\
\end{array}
\right. 
$$
and this is bounded at 0.\\
 Letting 
 $$B_3=-(m-1)\int_1^0[R^{-1}(s)s(1-f_d^2)\left(\begin{array}{c}
0\\
a_0 \\
0\\
b_0
\end{array}
 \right)]_3ds \hbox{  and  } \hat{X}_3=X_3\int_{0}^r[R^{-1}(s)s(1-f_d^2)\left(\begin{array}{c}
0\\
a_0 \\
0\\
b_0
\end{array}
 \right)]_3ds,$$
 we find 
 $$A_3(r)X_3=(A_3+B_3)X_3-(m-1)\hat{X}_3$$
with 
$$
\hat{X}_3=X_3O(1) \hbox{  at $+\infty$  }$$
and $\hat{X}_3$ is bounded at 0.\\ 
 For the fourth term,  
 $$[R^{-1}(s)s(1-f_d^2)\left(\begin{array}{c}
0\\
a_0 \\
0\\
b_0
\end{array}
 \right)]_4\sim_{+\infty}{-1\over 16n\sqrt2}{4d^2\sqrt2\over s}s^{-n}(-a_0+b_0)
 $$
 and 
$$
 \sim_0\left\{\begin{array}{c}
 s(O(s^{\gamma_1}a_0)+O(s^{2\gamma_1+\gamma_2}b_0))\hbox{  if  }(d,\gamma_1,\gamma_2)\in\mathcal{D}_1\hbox{  and  }\gamma_1+\gamma_2-2d-2<0\\
  s(O(s^{\gamma_1}a_0)+O(s^{\-\gamma_1+2d+2}b_0))\hbox{  if  }(d,\gamma_1,\gamma_2)\in\mathcal{D}_1\hbox{  and  }\gamma_1+\gamma_2-2d-2>0\\
 s\tau(s)(O(s^{\gamma_1})a_0+O(s^{\gamma_1+2d+2})b_0)\hbox{  if  }(d,\gamma_1,\gamma_2)\in\mathcal{D}_2
 \end{array}
 \right.
 $$
  Letting 
 $$B_4=-(m-1)\int_1^{0}[R^{-1}(s)s(1-f_d^2)\left(\begin{array}{c}
0\\
a_0 \\
0\\
b_0
\end{array}
 \right)]_4ds \hbox{  and  }\hat{X}_4=X_4\int_{0}^r[R^{-1}(s)s(1-f_d^2)\left(\begin{array}{c}
0\\
a_0 \\
0\\
b_0
\end{array}
 \right)]_4ds$$
 we find 
 $$A_4(r)X_4=(A_4+B_4)X_4-(m-1)\hat{X}_4.$$
 Then $\hat{X}_4=X_4O(1)$ at $+\infty$
and $\hat{X}_4$ is bounded at 0.\\
 Now, summing the four terms, and letting $C_i=A_i+B_i$, we find
  (\ref{eq:X0-1}) and (\ref{eq:X0-2}).\\
Since $X_0$ is bounded at 0, we have $C_2=C_4=0$. \\
But $X_0$ is bounded at $+\infty$ and $\hat{X}_i$ is bounded at $+\infty$, $i=1,2,3$. Since we have also $a_1>>\hat{a}_4$  at $+\infty$, we infer that $C_1=0$ and that $\hat{X}_4$ must be bounded at $+\infty$. Returning to the definition of $\hat{X}_4$, we must have
 $$\int_{0}^{+\infty}[R^{-1}(s)s(1-f_d^2)\left(\begin{array}{c}
0\\
a_0 \\
0\\
b_0
\end{array}
 \right)]_4ds=0,$$
therefore
 $$\hat{X}_4=X_4\int_{+\infty}^r[R^{-1}(s)s(1-f_d^2)\left(\begin{array}{c}
0\\
a_0 \\
0\\
b_0
\end{array}
 \right)]_4ds,$$
 that gives
 $$\hat{X}_4=X_4\int_{+\infty}^rs(1-f_d^2)[a_0\mathcal{C}_2+b_0\mathcal{C}_4]_4\sim_{+\infty}X_4\int_{+\infty}^r{-8\sqrt2\over-16n\sqrt2} s^{-n_0-n}{d^2\over s}ds.$$
Thus,
 \begin{equation}\label{eq:nveauX4}
\hbox{at $+\infty$ }\quad \hat{a}_4=a_4{-1\over 16n\sqrt2}{8d^2\sqrt 2\over n+n_0}r^{-n}a_0+o(r^{-n_0}).
 \end{equation}
  Since we have now
 $$X_0=C_3X_3-(m-1)\sum_{i=1}^4\hat{X}_i$$
 and since  $\hat{a}_1=O(r^{-n_0-4})$ and $\hat{a}_2=O(r^{-n_0-4})$, then $\hat{a}_1=o(a_0)$ and $\hat{a}_2=o(a_0)$  at $+\infty$. Consequently
\begin{equation}\label{eq:equiX0}
a_0+(m-1)\hat{a}_4\sim_{+\infty}C_3a_3-(m-1)\hat{a}_3
\end{equation}
Recalling (\ref{eq:nveauX4}) and recalling $n<n_0$, this implies that
 $$C_3-(m-1)\int_{0}^{+\infty}[R^{-1}(s)s(1-f_d^2)\left(\begin{array}{c}
0\\
a_0 \\
0\\
b_0
\end{array}
 \right)]_3ds=0$$
 and then
 $$C_3a_3-(m-1)\hat{a}_3=-(m-1)a_3\int_{+\infty}^{r}[R^{-1}(s)s(1-f_d^2)\left(\begin{array}{c}
0\\
a_0 \\
0\\
b_0
\end{array}
 \right)]_3ds.$$
Using (\ref{eq:3emetermeenlinfini}), we get
 \begin{equation}\label{eq:nveauX3}
\hbox{at $+\infty$}\quad C_3a_3-(m-1)\hat{a}_3=-(m-1)a_3{-1\over16\sqrt2}{-8d^2\sqrt2\over n( n-n_0)}r^{n-n_0}+o(r^{-n_0}).
 \end{equation}
 Finally, we sum (\ref{eq:nveauX4}) and (\ref{eq:nveauX3}) to get, by (\ref{eq:equiX0})
$$a_0(1+(m-1){-8d^2\over 16 n}{1\over n+n_0})\sim_{+\infty}-(m-1){8d^2\over 16 n}({1\over n-n_0})r^{-n_0}$$
 and thus
 $$(m-1){8d^2\over16 n}({-1\over n-n_0}+{1\over n+n_0})=1.$$
 But we have by (\ref{eq:n0})
 $$n^2_0-n^2=(-m+1)d^2.$$
After simplification by $m-1$, we get $n_0=n$, that gives $m=1$, that is in contradiction with the hypothesis $m<1$. So we deduce that $m=1$.\\
  The proof of (\ref{eq:X0-1}) and (\ref{eq:X0-2}) for $(d,\gamma_1,\gamma_2)\in\mathcal{D}_1$ and $\gamma_1+\gamma_2-2d-2=0$ is left to the reader. \\
  {\bf Proof of Corollary \ref{solutionbornee}.}\\
  By Theorem \ref{lierlespb} (ii), if there exists a bounded solution $\omega$, then there exists some eigenvalue tending to 1. So $m=1$. It remains to prove that $\omega_0=c\omega$, for some $c\neq 0$. But $\omega$ cannot have the least behavior at 0, otherwise it would blow up exponentially at $+\infty$. So, there exists $c\neq0$ such that $\omega\sim_0 c\omega_0$. If $\omega\neq c\omega_0$, then $\omega-c\omega_0$ has the least behavior at 0, and consequently blows up exponentially at $+\infty$. This cannot be true, because $\omega$ is bounded at $+\infty$ and, since $a_0\geq b_0\geq0$, the possible blowing up behavior at $+\infty$ for $\omega_0$ can only be polynomial. We can conclude that $\omega=c\omega_0$.
\section{The case $n\geq d+1$ : the proof of Theorem \ref{ngeqd+1}.}
Let
$\omega_1=(a_1,b_1)$ be the solution defined in Theorem \ref{solutionscontinuesen0} and $\eta_2=(u_2,v_2)$ be the solution defined in Theorem \ref{solutions continues en infty}. According to Theorem \ref{lesextremes}, $\omega_1\sim_{+\infty}(J_+,J_+)$ and $\eta_2$ has the greater blowing up behavior at 0. Let $\eta_3$ and $\eta_4$  be defined in Theorem \ref{solutions continues en infty} and having the intermediate behaviors at $+\infty$.
Let $\omega_3=(a_3,b_3)$ be defined in Theorem \ref{solutionscontinuesen0}. 
With these definitions, we can write
 $$\omega_3=C_1(n,d)\omega_1+C_2(n,d)\eta_2+C_3(n,d) \eta_3+C_4(n,d)\eta_4.$$
 Let us remark that $\omega_1$ and $\omega_3-C_1(n,d)\omega_1$ form a base of the bounded solutions at 0, and that $\omega_3-C_1(n,d)\omega_1=o(\omega_1)$  at $+\infty$. So the problem of the existence of some bounded solution is reduced to the problem $C_3(n,d)=0$.\\
 Supposing that there exists a bounded solution for $(n_0,d_0)$, $d_0>1$, $n_0\geq d_0+1$, we have, by Theorem \ref{essentiel}, $n_0\leq 2d_0-1$. From now on, $(n,d)$ is such that $1\leq d\leq d_0+1$ and $d\leq n\leq 2d$. Clearly, $(d,\vert n-d\vert, n+d)$ stays in a compact subset of $\mathcal{D}$. This is sufficient for the solutions $\eta_3$ and $\eta_4$ to be defined without ambiguity. The real numbers $C_i(n,d)$ defined above can be computed by means of determinants involving the four components $(a,a',b,b')(r)$ of the five solutions present, for a given $r>0$. Thus, $C_i$ is continuous wrt $(d,\gamma_1,\gamma_2)$ and consequently is continuous wrt $(d,n)$. $C_i$ is also differentiable wrt $\gamma_1$ and wrt $\gamma_2$ and therefore wrt $n$, since $n\geq d$.\\
  \begin{lemma}\label{deriverC3engamman}
With the notation above, if there exists $(n_0,d_0)$, $d_0\geq1$, $n_0\geq d_0+1$ such that $C_3(n_0,d_0)=0$, then  there exists a continuous map $d\mapsto n(d)$, defined for  $d<d_0$, closed to $d_0$ and verifying  $C_3(n(d),d)=0$.
\end{lemma}
 {\bf Proof  }  Let us prove that ${\partial C_3\over\partial{n}}(n_0,d_0)\neq0.$ If  ${\partial C_3\over\partial n}(n_0,d_0)=0,$ then 
 ${\partial\over\partial{n}}(\omega_3-C_1(n,d)\omega_1)(n_0,d_0)$ is bounded at  $+\infty$.
 Let us denote  $(a,b)=\omega_3-C_1(n,d)\omega_1$. Then $(a,b)$ verifies the system (\ref{eq:GL}), with $(n_0,d_0)$ in place of $(n,d)$, and $({\partial a\over\partial{n}}, {\partial b\over\partial{n}})(n_0,d_0)$ verifies also a system, obtained by differentiation wrt $n$, at $(n_0,d_0)$, that is 
 \begin{equation}\label{eq:systdiff}
\left\{\begin{array}{rl}
{\partial a\over\partial{n}}''+{1\over r}{\partial a\over\partial{n}}'-{(n-d)^2\over r^2}{\partial a\over\partial{n}}-2{(n-d)\over r^2}a-f_d^2{\partial b\over\partial{n}}&=-(1-2f_d^2){\partial a\over\partial{n}}\\
{\partial b\over\partial{n}}''+{1\over r}{\partial b\over\partial{n}}'-{(n+d)^2\over r^2}{\partial b\over\partial{n}}-2{(n+d)\over r^2}b-f_d^2{\partial a\over\partial{n}}&=-(1-2f_d^2){\partial b\over\partial{n}}
\end{array}
\right. 
\end{equation} 
  By combining the systems (\ref{eq:GL}) and (\ref{eq:systdiff}), for $(n_0,d_0)$, an integration by parts gives
$$\int_0^{+\infty}-2{n_0-d_0\over r}a^2-2{n_0+d_0\over r}b^2dr=0$$
 and we conclude that $a=b=0$, that is false.\\
  So, we have proved that ${\partial C_3\over\partial n}(n_0,d_0)\neq 0$. The Implicit Functions Theorem gives a continuous map  $d\mapsto n(d)$  such that $C_3(n(d),d)=0$, and defined in a neighborhood of $d_0$, with values in a neighborhood of $n_0$.\\
 {\bf The proof of Theorem \ref{ngeqd+1} completed.}\\
With the definitions given above,
let us define the set
 $$\mathcal{E}=\{d\geq1;\quad d\leq d_0+1;\quad \exists n\geq d+{1\over2},  \quad C_3(n,d)=0\}.$$
 If $d\in \mathcal{E}$, then $n\leq 2d-1$, by Theorem \ref{essentiel}. Thus, $\mathcal{E}$ is a closed subset of $[1,+\infty[$, thanks to the continuity of $C_3$ wrt $(n,d)$. 
Since $d_0\in\mathcal{E}$, $\mathcal{E}\neq\emptyset$ and we let $d_1=\inf\mathcal{E}$. Given that $d_1\in{\mathcal{E}}$, there exists $n_1\geq d_1+{1\over2}$ such that $C_3(n_1,d_1)=0$. According to Theorem \ref{n<d+1}, $n_1\geq d_1+1$. If $d_1>1$, we deduce from Lemma \ref{deriverC3engamman} that there exists $d<d_1$, sufficiently closed to $d_1$ in order to have $n(d)>d_1+{1\over2}$. Therefore $n(d)\geq d+{1\over2}$, which is in contradiction with $d_1=\inf\mathcal{E}$. This proves that $d_1=1$. But $1\not\in\mathcal{E}$, by Theorem \ref{essentiel}. 
  This contradiction proves the non existence of $(n_0,d_0)$ such that $n_0\geq d_0+1$ and $C_3(n_0,d_0)=0$.\\
  The proof of Theorem \ref{ngeqd+1} is complete.
\section{The proof of Theorem \ref{lierlespb} (i) and  of Theorem \ref{valeurproprenegative}.}  
{\bf Proof of Theorem \ref{lierlespb} (i).}
Let us define $n_0=\sqrt {{\gamma_2^2+\gamma_1^2\over 2}-\mu d^2}$.\\
     Let  $\omega_{\ep}=(a_{\ep},b_{\ep})\in\mathcal{H}_{\gamma_1}$  be an eigenvector associated to $\mu(\ep)$.  Using (\ref{eq:GLvp}), we write
 $${\mu(\ep)\over\ep^2}\int_0^1r(1-f^2)(a^2_{\ep}+b^2_{\ep})dr=\int_0^{1}(ra'^2_{\ep}+rb'^2_{\ep}+{\gamma^2_1\over r} \ a^2_{\ep}+{\gamma^2_2\over r} \ b^2_{\ep}+{r\over\ep^2}f^2( a_{\ep}+ b_{\ep})^2)dr.$$
 We use the definition (\ref{eq:defm0}) of $m_0(\ep)$ to get
 $${\mu(\ep)\over\ep^2}\int_0^1r(1-f^2)(a^2_{\ep}+b^2_{\ep})dr$$
 $$\geq {m_0(\ep)\over\ep^2}\int_0^1r(1-f^2)(a^2_{\ep}+b^2_{\ep})dr+\int_0^1({\gamma^2_1-d^2\over r}  a^2_{\ep}+{\gamma^2_2-d^2\over r}  b^2_{\ep}
 +{r\over \ep^2}{f}^2( a_{\ep}+b_{\ep})^2)dr.$$
 Now, we use the trick of TC Lin (see \cite{TCL1}). Letting $\tilde b_{\ep}=\tau \tilde a_{\ep}$, we consider the map
 \begin{equation}\label{eq:defH}
 H:\tau\mapsto {\gamma^2_1-d^2\over r} +{\gamma^2_2-d^2\over r}\tau^2+r{f_d} ^2(1+\tau)^2
 \end{equation}
 and we minimize this map. The minimum is attained for $\tau_0$ verifying
 $$\tau_0({\gamma^2_2-d^2\over r}+rf^2_d)+rf^2_d=0\quad\hbox{  
 and   }\quad
 1+\tau_0={{\gamma^2_2-d^2\over r}/({\gamma^2_2-d^2\over r}+rf^2_d})$$
 and consequently
 $$H(\tau_0)={\gamma^2_1-d^2\over r}+({rf^2_d\over{\gamma^2_2-d^2\over r}+rf^2_d})^2({\gamma^2_2-d^2\over r})
+rf_d^2({{\gamma^2_2-d^2\over r}
   \over {\gamma^2_2-d^2\over r}+rf_d^2})^2.$$
We have
$$H(\tau_0)\sim_{r\rightarrow+\infty}({\gamma^2_1+\gamma^2_2-2d^2)/ r}.\quad\hbox{
Moreover, for all $\tau>0,$  }\quad H(\tau)\geq H(\tau_0).$$
Since ${\gamma_1^2+\gamma_2^2\over2}-d^2>0,$  there exists some constants $C_1>0$ and $R_0>0$, independent of $\tau$, such that for all $\tau>0$
 $$H(\tau)\geq {C_1\over r}\hbox{  for all $r>R_0$}.$$
 Then, for all  $R>R_0$ and all  $\ep<{1\over R}$, we write
 $$\int_0^{1\over\ep}H(r)\tilde{a}^2_{\ep}(r)dr\geq\int_0^{R_0}H(r)\tilde{a}^2_{\ep}(r)dr+\int_{R_0}^RH(r)\tilde{a}^2_{\ep}(r)dr.$$
 Now $a_0$ blows up exponentially at $+\infty$, or as $r^{n_0}$. We can choose $R_0$ large enough and a constant $C_2>0$ to have also
 $$a_0^2(r)\geq C_2 ({e^{\sqrt2 r}\over\sqrt r})^2
 \hbox{  or  }
C_2 r^{2n_0}
 \hbox{  for all $r>R_0$}.$$
 Since $\tilde a_{\ep}\rightarrow a_0$ as $\ep\rightarrow0$, uniformly in $[0,R_0]$, we can chose $\ep_0$ such that for all $\ep<\ep_0$
 $$\int_0^{R_0}H(r)\tilde{a}^2_{\ep}(r)dr\geq {1\over2}\int_0^{R_0}H(r){a}^2_{0}(r)dr.$$
Moreover, for all $R>R_0$, $\tilde a_{\ep}\rightarrow a_0$ as $\ep\rightarrow0$, uniformly in $[R_0,R]$.
 Then, there exists $\ep(R)$ such that for all $\ep<\ep(R)$ we have
 $$\int_{R_0}^RH(r)\tilde{a}^2_{\ep}(r)dr\geq {C_2\over2}\int_{R_0}^{R}{1\over r} r^{2n_0} dr\quad\hbox{   or    }\quad
 \int_{R_0}^RH(r)\tilde{a}^2_{\ep}(r)dr\geq {C_2\over2}\int_{R_0}^{R}{1\over r} ({e^{\sqrt2 r}\over\sqrt r})^2 dr.$$
 And finally, for $\ep<\ep(R)$,  we have
 $$({\mu(\ep)-m_0(\ep)\over\ep^2})\int_0^1r(1-f^2)(a^2_{\ep}+b^2_{\ep})dr\geq{1\over2}\int_0^{R_0}H(r){a}^2_{0}(r)dr+$$
$$+\left\{\begin{array}{c}{C_1C_2\over2}\int_{R_0}^{R}{1\over r} r^{2n_0},\quad\hbox{if $(a_0,b_0)\sim_{+\infty}(r^{n_0},-r^{n_0})$}\\
{C_1C_2\over2}\int_{R_0}^{R}{1\over r}({e^{\sqrt2 r}\over\sqrt r})^2dr, \quad\hbox{if $(a_0,b_0)\sim_{+\infty}(J^{+},J^{+})$}
\end{array}
\right.
$$
 where $C_1$ and $C_2$, given above, are independent of $R$ and $\ep$. But we can choose $R$ such that the lhs is positive.\\
We deduce that
 ${\mu(\ep)-m_0(\ep})>0.$
 Then we use Theorem \ref{lambda0} (i), that gives ${m_0(\ep)-1\over\ep^2}\geq C$ and consequently ${\mu(\ep)-1\over\ep^2}\geq C$.
 The lemma is proved.\\
 {\bf The proof of Theorem \ref{valeurproprenegative}.}
The proof for $n=2$ and $d=2$ is originally  in  \cite{Miro}. \\
  For real numbers $d\geq1$ and $n\geq1$, let $x={f'_d\over r^{n-1}}$ and $y=d{f_d\over r^n}$. A calculation gives
 \begin{equation}
 \left\{\begin{array}{rl}-(rx')'+{\gamma^2\over r}x-{\xi^2\over r}y-r(1-3f_d^2)x&=-2{n-1\over r^{n-1}}f_d(1-f_d^2)\\
 -(ry')'+{\gamma^2\over r}y-{\xi^2\over r}x-r(1-f_d^2)y&=0
 \end{array}
 \right.
 \end{equation}
 For $a={x+y\over2}$ and $b={x-y\over2}$, we deduce that
 \begin{equation}\label{eq:syst}
 \left\{\begin{array}{rl}-(ra')'+{\gamma_1^2\over r}a+f_d^2b-r(1-2f_d^2)a&=-{n-1\over r^{n-1}}f_d(1-f_d^2)\\
 -(rb')'+{\gamma_2^2\over r}b+f_d^2a-r(1-2f_d^2)b&=-{n-1\over r^{n-1}}f_d(1-f_d^2)
 \end{array}
 \right.
 \end{equation}
 where, as usual, $\gamma_1=\vert n-d\vert$, $\gamma_2=n+d$, $\gamma^2={\gamma_1^2+\gamma_2^2\over2}$ and $\xi^2={\gamma_2^2-\gamma_1^2\over2}$.\\
 We verify that
 $$x\sim_0 y\sim_0 dr^{d-n}+O(r^{d-n+2}) \hbox{  and, at $+\infty$, }x=O(r^{-n}),\quad y=O(r^{-n}),$$
 and consequently that
 $$a\sim_02dr^{d-n}+O(r^{d-n+2})\hbox{  and  }b\sim_0 O(r^{d-n+2}).$$
 let us suppose that $d\geq1$ and that $1< n<d+1$. 
We can multiply the system (\ref{eq:syst}) and integrate by parts. We obtain that
 
 $$\int_0^{+\infty}(ra'^2+rb'^2+{\gamma^2_1\over r}a^2+{\gamma^2_2\over r}b^2+rf_d^2(a+b)^2-r(1-f_d^2)(a^2+b^2))dr$$
 $$=\int_0^{+\infty}{-(n-1)\over r^{n-1}}f_d(1-f_d^2)(a+b)dr
 $$
This gives
$$ {\int_0^{+\infty}(ra'^2+rb'^2+{\gamma^2_1\over r}a^2+{\gamma^2_2\over r}b^2+rf_d^2(a+b)^2)dr\over\int_0^{+\infty}r(1-f_d^2)(a^2+b^2)dr}
=1-C_n$$
with 
$$C_n={\int_0^{+\infty}{(n-1)\over r^{n-1}}f_d(1-f_d^2)(a+b)dr\over\int_0^{+\infty}r(1-f_d^2)(a^2+b^2)dr}>0.$$
 Now we use an approximation argument, valid as soon as $n>0$. For example for a given constant $0<N<1$ we define\\
 $({a}_{\ep},{b}_{\ep})(r)=\left\{\begin{array}{c}
(a,b)({r\over\ep})\hbox{ in $[0,N]$}\\
=(a(r){(1-r)^2\over (1-N)^2},b(r){(1-r)^2
\over (1-N)^2})\hbox{ in $[N,1]$}
\end{array}.
\right.$
We have that  $(a_{\ep},b_{\ep})\in\mathcal{H}_{\vert n-d\vert}$ and that
$$
{\int_0^{1}(ra_{\ep}'^2+rb_{\ep}'^2+{\gamma^2_1\over r}a_{\ep}^2+{\gamma^2_2\over r}b_{\ep}^2+r{1\over\ep^2}f^2(a_{\ep}+b_{\ep})^2)dr\over{1\over\ep^2}\int_0^{1}r(1-f^2)(a_{\ep}^2+b_{\ep}^2)dr}$$
$$= {\int_0^{N\over\ep}(ra'^2+rb'^2+{\gamma^2_1\over r}a^2+{\gamma^2_2\over r}b^2+rf_d^2(a+b)^2)dr+O(\ep^{2n})\over\int_0^{N\over\ep}r(1-f_d^2)(a^2+b^2)dr+O(\ep^{2n})}
\rightarrow 1-C_n,\hbox{  as $\ep$ tends to 0}.
$$
We deduce that, if $1<n<d+1$,
$m_{d-n,d+n}(\ep)<1-{C_n\over2}$,
for $\ep$ small enough and the proof of Theorem \ref{valeurproprenegative} is complete.

 \end{document}